\documentclass[11pt]{article}
\usepackage{amssymb,amsmath,latexsym}
\begin{document}

\begin{doublespace}

\newtheorem{thm}{Theorem}[section]
\newtheorem{lemma}[thm]{Lemma}
\newtheorem{defn}[thm]{Definition}
\newtheorem{prop}[thm]{Proposition}
\newtheorem{corollary}[thm]{Corollary}
\newtheorem{remark}[thm]{Remark}
\newtheorem{example}[thm]{Example}
\numberwithin{equation}{section}
\def\ee{\varepsilon}
\def\qed{{\hfill $\Box$ \bigskip}}
\def\NN{{\cal N}}
\def\AA{{\cal A}}
\def\MM{{\cal M}}
\def\BB{{\cal B}}
\def\CC{{\cal C}}
\def\LL{{\cal L}}
\def\DD{{\cal D}}
\def\FF{{\cal F}}
\def\EE{{\cal E}}
\def\QQ{{\cal Q}}
\def\RR{{\mathbb R}}
\def\R{{\mathbb R}}
\def\L{{\bf L}}
\def\E{{\mathbb E}}
\def\F{{\bf F}}
\def\P{{\mathbb P}}
\def\N{{\mathbb N}}
\def\eps{\varepsilon}
\def\wh{\widehat}
\def\wt{\widetilde}
\def\pf{\noindent{\bf Proof:} }
\def\pfof{\noindent{\bf Proof} }

\def\sA {{\cal A}} \def\sB {{\cal B}} \def\sC {{\cal C}}
\def\sD {{\cal D}} \def\sE {{\cal E}} \def\sF {{\cal F}}
\def\sG {{\cal G}} \def\sH {{\cal H}} \def\sI {{\cal I}}
\def\sJ {{\cal J}} \def\sK {{\cal K}} \def\sL {{\cal L}}
\def\sM {{\cal M}} \def\sN {{\cal N}} \def\sO {{\cal O}}
\def\sP {{\cal P}} \def\sQ {{\cal Q}} \def\sR {{\cal R}}
\def\sS {{\cal S}} \def\sT {{\cal T}} \def\sU {{\cal U}}
\def\sV {{\cal V}} \def\sW {{\cal W}} \def\sX {{\cal X}}
\def\sY {{\cal Y}} \def\sZ {{\cal Z}}

\def\bA {{\mathbb A}} \def\bB {{\mathbb B}} \def\bC {{\mathbb C}}
\def\bD {{\mathbb D}} \def\bE {{\mathbb E}} \def\bF {{\mathbb F}}
\def\bG {{\mathbb G}} \def\bH {{\mathbb H}} \def\bI {{\mathbb I}}
\def\bJ {{\mathbb J}} \def\bK {{\mathbb K}} \def\bL {{\mathbb L}}
\def\bM {{\mathbb M}} \def\bN {{\mathbb N}} \def\bO {{\mathbb O}}
\def\bP {{\mathbb P}} \def\bQ {{\mathbb Q}} \def\bR {{\mathbb R}}
\def\bS {{\mathbb S}} \def\bT {{\mathbb T}} \def\bU {{\mathbb U}}
\def\bV {{\mathbb V}} \def\bW {{\mathbb W}} \def\bX {{\mathbb X}}
\def\bY {{\mathbb Y}} \def\bZ {{\mathbb Z}}
\def\R {{\mathbb R}} \def\RR {{\mathbb R}}
\def\n{{\bf n}}

\title{Minimal thinness for subordinate Brownian motion in half-space}

\author{{\bf Panki Kim}\thanks{This work was supported by Basic Science Research Program through the National Research Foundation of Korea(NRF) funded by the Ministry of Education, Science and Technology(0409-20110087)} \quad {\bf Renming Song} \quad and
\quad {\bf Zoran Vondra\v{c}ek}\thanks{Supported in part by the MZOS
grant 037-0372790-2801.} }

\date{July 26, 2011}

\maketitle

\begin{abstract}
We study minimal thinness in the half-space $H:=\{x=(\wt{x}, x_d):\,
\wt{x}\in \R^{d-1}, x_d>0\}$
for a large class of
rotationally invariant L\'evy processes, including symmetric stable
processes and sums of Brownian motion and independent stable
processes. We show that the same test for the minimal thinness of a
subset of $H$ below the graph of a nonnegative
Lipschitz function is valid for all processes
in the considered class. In the classical case of Brownian motion
this test was proved by Burdzy.
\end{abstract}

\noindent {\bf AMS 2010 Mathematics Subject Classification}: Primary
60J50, 31C40; Secondary 31C35, 60J45, 60J75.

\noindent {\bf Keywords and phrases:} Minimal thinness, subordinate
Brownian motion, boundary Harnack principle, Green function, Martin
kernel

\section{Introduction}
Minimal thinness is a notion that describes the smallness of a set
at a boundary point. More precisely, let $D$ be a domain in $\R^d$,
$d\ge 2$, let  $\partial^M D$ (respectively $\partial^m D$) denote
its Martin boundary (respectively minimal Martin boundary) with
respect to Brownian motion, and let $M^D(x,z)$, $x\in D$, $z\in
\partial^M D$, be the corresponding  Martin kernel with respect to
Brownian motion. For $A\subset D$, let $\wh{R}^A_{M^D(\cdot, z)}$
denote the balayage of $M^D(\cdot, z)$ onto $A$.
The set $A$ is said to be  minimally thin in $D$ at $z\in
\partial^m D$ with respect to Brownian motion if $\wh{R}^A_{M^D(\cdot,
 z)}\neq M^D(\cdot, z)$. The concept of minimal thinness in the context of
classical potential theory was introduced and studied by Na\"{\i}m
in \cite{Nai}; for a recent exposition see \cite[Chapter 9]{AG}.
A probabilistic interpretation of minimal thinness is due to Doob,
see, e.g., \cite{Doo}: $A\subset D$ is minimally thin in $D$
 with respect to Brownian motion at $z\in \partial^m D$ if there
exists a point $x\in D$ such that with positive probability the
$M^D(\cdot,z)$-conditioned Brownian motion starting from $x$ does not
hit $A$.

We recall now two results about minimal thinness in the half-space
$H:=\{x=(\wt{x}, x_d):\, \wt{x}\in \R^{d-1}, x_d>0\}$, $d\ge 2$,
with respect to Brownian motion. The Martin boundary of $H$ with
respect to Brownian motion  can be identified with $\partial H \cup
\{\infty\}$, where $\partial H=\{(\wt{x}, 0):\,  \wt{x}\in
\R^{d-1}\}$, and all boundary points are minimal. The first result
is due to Beurling \cite{Beu} in the case $d=2$ and Dahlberg
\cite{Dal} in the case $d\ge 3$. By $B(z,r)$ we denote the ball
centered at $z\in \R^d$ with radius $r>0$.

\begin{thm}\label{t:dahlberg}
Let $A$ be a Borel subset of $H$ and assume that
\begin{equation}\label{c:criterion1}
\int_{A\cap B(0,1)} |x|^{-d}\, dx =\infty\, .
\end{equation}
Then $A$ is not minimally thin in $H$ with respect to Brownian
motion at $z=0$.
\end{thm}

The second result is a test for the minimal thinness of a subset of $H$
below the graph of a nonnegative Lipschitz function originally
proved by Burdzy \cite{Bur} using a probabilistic approach. An
alternative proof using Theorem \ref{t:dahlberg} was given by
Gardiner \cite{Gar}.

\begin{thm}\label{t:burdzy}
Let $f:\R^{d-1}\to [0,\infty)$ be a Lipschitz function with
Lipschitz constant $a>0$. The set $A:=\{x=(\wt{x},x_d)\in H:\,
0<x_d\le f(\wt{x})\}$ is minimally thin in $H$ with respect to
Brownian motion at $z=0$ if and only if
\begin{equation}\label{c:criterion2}
\int_{\{|\wt{x}|<1\}}f(\wt{x})|\wt{x}|^{-d}\, d\wt{x} <\infty\, .
\end{equation}
\end{thm}

The goal of this paper is to show that the above two theorems are
still valid in exactly the same form when Brownian motion is
replaced with a wide class of rotationally invariant L\'evy processes
(see Theorems \ref{t:main1}--\ref{t:main2}). The precise
description of this class will be given in the next section -- for
now it suffices to know that it includes rotationally invariant
$\alpha$-stable processes, $\alpha\in (0,2)$. The Martin boundary
theory for Hunt processes admitting a dual process (and satisfying
an additional hypothesis) was developed by Kunita and Watanabe
\cite{KW}, while the concept of minimal thinness for such processes
was studied by F\"ollmer \cite{Fol}. To the best of our knowledge no
concrete criteria
for minimal thinness in the spirit of
Theorems \ref{t:dahlberg}--\ref{t:burdzy} have been obtained for
 any discontinuous processes, not even the symmetric stable ones.
Time is now ripe for such results due to the recent progress in the
potential theory of rotationally invariant L\'evy processes, in
particular subordinate Brownian motions.
Our proofs of the analogs of Theorems
\ref{t:dahlberg}--\ref{t:burdzy} will heavily rely on the very
recent work
\cite{KSV1, KSV3, KSV2, KSV4}
where a boundary Harnack principle and sharp estimates of the Green
function of certain subordinate Brownian motions were obtained.

We find the conclusion of our main result, Theorem \ref{t:main2},
surprising since the test \eqref{c:criterion2} is the same for all
processes in the considered class. In particular, for symmetric
$\alpha$-stable processes, the criterion for the minimal thinness of
the set $A$ in $H$ at $z=0$ does not depend on the index of
stability $\alpha$. This is in contrast with the following criterion
for the thinness of thorns: Let $f:[0,\infty)\to [0,\infty)$ be an
increasing function such that $f(r)>f(0)$ for all $r>0$, and
$f(r)/r$ is non-decreasing for sufficiently small $r>0$. Let
$A:=\{x\in H:\, |\wt{x}| < f(x_d)\}$. Then $A$ is
thin in $H$ at $0$ with
respect to Brownian motion if and only if
\begin{eqnarray*}
& &\int_0^1 \left(\frac{f(r)}{r}\right)^{d-3}\, \frac{dr}{r}<\infty\, ,
\qquad d\ge 4\, ,\\
& &\int_0^1 \Big| \log\frac{f(r)}{r}\Big|^{-1}\frac{dr}{r} <\infty\,
,\qquad d=3\, .
\end{eqnarray*}
On the other hand, for $\alpha\in (0,2)$, $A$ is thin at $0$ with
respect to the symmetric $\alpha$-stable process if and only if
$$
\int_0^1 \left(\frac{f(r)}{r}\right)^{d-\alpha-1}\,
\frac{dr}{r}<\infty,\qquad d\ge 3.
$$
The above criterion for the thinness of thorns
 in $H$ at $0$ with
respect to
$\alpha$-stable processes can be proved by using Wiener's test for
stable processes \cite[Corollary 4.17]{BH} and by slightly modifying
the proofs in \cite[pp. 67--69]{PS} (by changing cylindrical
surfaces to full cylinders).

This paper is organized as follows. In the next section we precisely
describe the class of subordinate Brownian motions for which we
will study minimal thinness in the half-space, recall from
\cite{KSV3, KSV2, KSV4} relevant results on the boundary Harnack
principle for those processes, and derive necessary estimates for
the Green function of the half-space (for points close to the
boundary and to each other). In Section 3 we use these Green
function estimates to obtain two-sided estimates on the Martin
kernel for points close to the origin and to each other. These
estimates and the boundary Harnack principle suffice to identify the
finite part of the minimal Martin boundary of $H$ with $\partial H$.
This result may be of independent interest. We then show how the
studied processes fit in the framework of minimal thinness in
\cite{Fol} and recall both the potential-theoretic and the
probabilistic definitions of minimal thinness. In the last section
we state and prove Theorems \ref{t:main1}--\ref{t:main2}, the
analogs of Theorems \ref{t:dahlberg}--\ref{t:burdzy} for our
processes. The main ingredients of the proof are Lemma
\ref{l:sjogren} and Proposition \ref{p:sjogren} which are
generalizations of \cite[Theorem 2]{Sj2}.
 Instead of the estimates of the classical Green function and Martin
kernel we use our estimates from Sections 2 and 3 for points close
to the origin and the boundary Harnack principle for points away
from the boundary.

We will use the following conventions in this paper. The values of
the constants $C_{1}(R), \dots, C_5(R)$, depending only on
$d$, $R>0$ and the Laplace exponent of the subordinator, will remain the same throughout this
paper, while the constants $c, c_0, c_1, c_2, \dots$  stand for
constants whose values are unimportant and which may change from one
appearance to another. All constants are positive finite numbers.
We assume  $d \ge 2$ and the dependence of the constants on the
dimension $d$ may not be mentioned explicitly.

For two nonnegative functions $f, g$, $f(t) \sim g(t)$, $t \to 0$
($f(t) \sim g(t)$, $t \to \infty$, respectively) means that $ \lim_{t \to
0} f(t)/g(t) = 1$ ($\lim_{t \to \infty} f(t)/g(t) = 1$,
respectively). On the other hand, $f(t) \asymp g(t)$, $t \to 0$
($f(t) \asymp g(t)$, $t \to \infty$, respectively) means that the
quotient  $f(t)/g(t)$ stays bounded between two positive
constants as $t \to 0$ (as $t \to \infty$, respectively).
Simply, $f\asymp g$ means that the
quotient  $f(t)/g(t)$ stays bounded between two positive
constants on their
common domain of definitions.

 For any open set $U$, we denote by $\delta_U (x)$ the
distance between $x$ and the complement of $U$, i.e.,
$\delta_U(x)=\text{dist} (x, U^c)$.  We will use $dx$ to denote the
Lebesgue measure in $\R^d$. For a Borel set $A\subset \R^d$, we also
use $|A|$ to denote its Lebesgue measure and $\mathrm{diam}(A)$ to
denote the diameter of the set $A$.
Finally, we will use ``$:=$" to denote a definition, which is
read as ``is defined to be".

\section{Preliminaries on subordinate Brownian motion}

In this section we will first describe a class of subordinate
Brownian motions and their potential theory. Recall that a
subordinator $S=(S_t)_{t\ge 0}$ is simply a nonnegative L\'evy
process
with $S_0=0$.
 The Laplace exponent of $S$ is a function
$\phi:(0,\infty)\to (0,\infty)$ having the representation
\begin{equation}\label{e:bernstein}
\phi(\lambda)= a\lambda +\int_{(0,\infty)} (1-e^{-\lambda t})\, \eta(dt)\, ,
\end{equation}
where $a\ge 0$ is the drift and $\eta$ the L\'evy measure of $S$,
i.e., a measure on $(0,\infty)$ satisfying
$\int_{(0,\infty)}(1\wedge t)\, \eta(dt)<\infty$. The Laplace
exponent $\phi$ determines the distribution of $S_t$ through the
formula $\E[\exp\{-\lambda S_t\}]=\exp\{-t\phi(\lambda)\}$. Formula
\eqref{e:bernstein} shows that $\phi$ is a Bernstein function,
i.e.~a nonnegative $C^{\infty}$ function on $(0,\infty)$ satisfying
$(-1)^{n-1} \phi^{(n)}\ge 0$ for all $n\ge 1$. Since the sample
paths of $S$ are nondecreasing functions, the subordinator $S$ can
serve as a stochastic time-change. More precisely, let
$Y=(Y_t,\P_x)_{t\ge 0, x\in \R^d}$ be a Brownian motion in $\R^d$
independent of $S$ with
$$
\E\left[e^{i\xi(Y_t-Y_0)}\right]=e^{-t{|\xi|^2}} \quad \xi\in \R^d,
t>0.
$$
The stochastic process $X=(X_t,\P_x)_{t\ge 0, x\in \R^d}$ defined by
the formula  $X_t:=Y_{S_t}$ is called a subordinate Brownian motion.
It is a rotationally invariant L\'evy process in $\R^d$ with
characteristic exponent $\Phi(\xi)=\phi(|\xi|^2)$ and infinitesimal
generator $-\phi(-\Delta)$. Here $\Delta$ denotes the Laplacian and
$\phi(-\Delta)$ is defined through functional calculus.

A Bernstein function $\phi$ is a complete Bernstein function if its
L\'evy measure $\eta$ has a completely monotone density, which will
be denoted by $\eta(t)$.
 We will
consider the following class of subordinate Brownian motions
determined mainly by the asymptotic behavior at infinity of the Laplace exponents of the corresponding subordinators:

\noindent {\bf Hypothesis (H)}: $d\ge 2$, $\phi$ is a complete Bernstein function, and
there exists $\alpha\in (0,2]$ such that $\phi(\lambda)\asymp \lambda^{\alpha/2} \ell(\lambda)$
as $\lambda \to \infty$, where $\ell:(0,\infty)\to (0,\infty)$ is measurable,
locally bounded above and below by positive constants,
and slowly varying at $\infty$.
Additionally,
\begin{itemize}
    \item in case $\alpha\in (0,2)$ and $d=2$, assume that there exists $\gamma <1$ such that
    $\liminf_{\lambda\to 0}\phi(\lambda)/\lambda^{\gamma}>0$;
    \item in case $\alpha=2$, assume that
    $d\ge 3$, $\phi$ has a positive drift $a$ and the L\'evy density $\eta$
    of $\phi$ satisfies the following condition: for any $K>0$, there exists
    $c=c(K)>1$ such that
    \begin{equation}\label{H:11}
    \eta(t)\le c\, \eta(2t), \qquad t\in (0, K).
    \end{equation}
    In this case one can take $\ell\equiv 1$.
\end{itemize}

\medskip

In case the above hypothesis holds true we will say that subordinate Brownian motion $X$ satisfies {\bf (H)}.
It is easy to check that in this case $X$ is transient.
Moreover, in this case the potential measure $U$ of
the corresponding subordinator $S$ has a density $u$ which is also
completely monotone (see, e.g., \cite[III, Theorem 5]{Ber}, \cite[Corollary 5.4 and Corollary 5.5]{BBKRSV} and \cite[Remark 10.6]{SSV}).

In the case $\alpha\in (0, 2)$, \eqref{H:11} is a consequence of the asymptotic behavior of $\phi$ at
infinity given in the first sentence of {\bf (H)}, see \cite[Theorem 2.10]{KSV3}.

Subordinate Brownian motions satisfying {\bf (H)} and with $\alpha \in (0,2)$  were studied in
\cite{KSV1, KSV3, KSV2}.
Such a subordinate Brownian motion $X$  is a purely
discontinuous L\'evy process in $\R^d$ with characteristic exponent
$\Phi$ satisfying $\Phi(\xi)\asymp |\xi|^{\alpha}\ell(|\xi|^2)$,
$|\xi|\to \infty$, $\xi\in \R^d$. This class of processes includes
$\alpha$-stable processes, corresponding to
$\phi(\lambda)=\lambda^{\alpha/2}$, relativistic $\alpha$-stable
processes, corresponding to
$\phi(\lambda)=(\lambda+m^{2/\alpha})^{\alpha/2}-m$, $m>0$, sums of
independent $\alpha$-stable and $\beta$-stable processes,
corresponding to
$\phi(\lambda)=\lambda^{\alpha/2}+\lambda^{\beta/2}$,
$0<\beta<\alpha<2$, and many others (see \cite{KSV2} for further
examples).

Subordinate Brownian motions satisfying {\bf (H)} with $\alpha=2$ and $d\ge 3$ were studied in
\cite{KSV4}. This class of processes includes independent sums of Brownian motion
and $\beta$-stable processes corresponding to
$\Phi(\xi)=a|\xi|^2+b^{\beta} |\xi|^{\beta}$ with $a, b>0$, and
independent sums of Brownian motion
and relativistic $\beta$-stable processes corresponding to $\Phi(\xi)=a|\xi|^2+(\lambda+m^{2/\beta})^{\beta/2}-m$ with $a, m>0$, and many others.

Let us first consider one-dimensional subordinate Brownian
motions. Suppose that $B=(B_t)_{t\ge 0}$ is a Brownian motion in
$\R$, independent of $S$, with
$$
\E\left[e^{i\theta(B_t-B_0)}\right]=e^{-t\theta^2}, \qquad\,
  \theta\in \R, \ t>0.
$$
The subordinate Brownian motion $Z=(Z_t)_{t\ge 0}$ in $\R$ defined
by $Z_t:=B_{S_t}$
 is a symmetric L\'evy process with characteristic
exponent $\Phi(\theta)=\phi(\theta^2)$, $\theta\in \R$. Define
$\overline{Z}_t:=\sup\{0\vee Z_s:0\le s\le t\}$ and let $L=(L_t:\,
t\ge 0)$ be a local time of $\overline{Z}-Z$ at $0$. $L$ is also
called a local time of the process $Z$ reflected at the supremum.
Then the right continuous inverse $L^{-1}_t$ of $L$ is a
 subordinator and is called the ladder time process
of $Z$. The process $\overline{Z}_{L^{-1}_t}$ is also a
 subordinator and is called the ladder height
process of $X$. (For basic properties of the ladder time and ladder
height processes, we refer the readers to \cite[Chapter 6]{Ber}.)
Let $\chi$ denote the Laplace exponent of the ladder height process
of $Z$, and let $V$ be its potential measure. By a slight abuse of
notation we also use $V$ to denote the function $V(t)=V((0,t))$,
$t>0$.

From now on we assume that the process $X$ is a subordinate Brownian
motion satisfying {\bf (H)}. Since $X$ is transient, it has a
Green function $G(x, y)$ given by
$$
G(x,y)=\int_0^{\infty}(4\pi t)^{-d/2}e^{-\frac{|x-y|^2}{4t}}u(t)\, dt\, ,
$$
where $u$ is the potential density of the subordinator $S$.
If we define
$$
G(r)=\int_0^{\infty}(4\pi t)^{-d/2}e^{-\frac{r^2}{4t}}u(t)\, dt\,
\quad r>0,
$$
then $G(\cdot)$ is a non-increasing function on $(0, \infty)$ and
$G(x,y)=G(|x-y|)$ for all $x, y\in \R^d$.

When $X$ is a subordinate Brownian motion satisfying {\bf (H)}, the Green function $G(x,y)$ of $X$
satisfies the following sharp estimates
\begin{equation}\label{e:asymp-G-I}
G(x,y)\asymp \frac{1}{|x-y|^d \phi(|x-y|^{-2})}, \qquad |x-y|\to 0.
\end{equation}
The Laplace exponent $\chi$ of the ladder height
process of $Z$ is a complete Bernstein function, the
function $V$ is a smooth function and satisfies
\begin{equation}\label{e:asymp-V-I}
V(t) \asymp \phi(t^{-2})^{-1/2}\, ,  \qquad t\to 0\, .
\end{equation}
For these two results see  \cite{KSV1, KSV3, KSV2} in case $\alpha\in (0,2)$,
and \cite{KSV4}
in case $\alpha=2$. In fact, when $\alpha=2$, we have more precisely,
\begin{eqnarray}
G(x,y)\sim \frac{\Gamma(d/2-1)}{4a\pi^{d/2}|x-y|^{d-2}}\, ,&& \qquad |x-y|\to 0\, ,\label{e:asymp-G-II}\\
V(t)\sim a\,t\, ,&& \qquad t\to 0\, .\label{e:asymp-V-II}
\end{eqnarray}

We record two consequences of estimates \eqref{e:asymp-G-I} and \eqref{e:asymp-V-I}.
\begin{prop}\label{p:GoverV}
Suppose that $X$ is a subordinate Brownian motion satisfying {\bf (H)}. Let $R>0$.

\noindent (i)
There exists a constant $C_1(R)=C_1(d,\phi, R)>1$
such that for all $x,y\in H$ satisfying $|x-y|<R$ it holds that
\begin{equation}\label{e:GoverV}
C_1(R)^{-1}|x-y|^{-d}\le \frac{G(|x-y|)}{V(|x-y|)^2}\le C_1(R)|x-y|^{-d}\, .
\end{equation}

\noindent
(ii)
There exists a constant $C_2(R)=C_2(d,\phi, R)>1$ such that
\begin{equation}\label{e:green-renewal}
C_2(R)^{-1}\le V(t)^{-2}\int_{B(0,t)}G(0,x)\, dx \le C_2(R)\, ,
\quad 0<t \le R\, .
\end{equation}
\end{prop}
\pf Fix $R>0$. We have by
\eqref{e:asymp-G-I} and \eqref{e:asymp-V-I} that for $|x-y|<R$
$$
G(|x-y|)\asymp \frac1{|x-y|^d \phi(|x-y|)^{-2})} \quad \text{and} \quad V(|x-y|)\asymp \phi(|x-y|^2)^{-1/2}.
$$
Now \eqref{e:GoverV} follows immediately with the constant depending only on $R, d, \phi$.

For part (ii) note that for $0<t\le R$,
$$
\int_{B(0,t)}G(0,x)\, dx=c_1 \int_0^t r^{d-1}G(r)\, dr \asymp
\int_0^t r^{d-1} \frac{1}{r^{d-\alpha}\ell(r^{-2})}\, dr
\asymp \frac{t^{\alpha}}{\ell(t^{-2})}\, ,
$$
where in case $\alpha\in (0,2)$ the last asymptotic equality follows by a property of slowly varying function $\ell$
(see \cite[Theorem 1.5.11]{BGT} - Karamata's theorem), while in the case $\alpha=2$ the last asymptotic equality is trivial.
Combining this with \eqref{e:asymp-V-I} we obtain \eqref{e:green-renewal}. \qed

A nonnegative function $h$ on $\R^d$ is harmonic in an open set
$D\subset \R^d$  with respect to
$X$ if for every open
set $B$ such that $B\subset \overline{B}\subset D$ and every $x\in
B$ it holds that $h(x)=\E_x[h(X_{\tau_B})]$. Here
$\tau_B=\inf\{t>0:\, X_t\notin B\}$ is the
first exit time of $X$ from
$B$. A nonnegative function $h$ on $\R^d$ is regular harmonic in an
open set $D\subset \R^d$  with respect to
$X$ if
$h(x)=\E_x[h(X_{\tau_D})]$ for all $x\in D$. The following Harnack principle is a consequence
of \cite[Theorem 4.7]{KSV3} when $\alpha\in (0,2)$, and
\cite[Theorem 4.5]{RSV} when $\alpha=2$.

\begin{thm}[Harnack inequality]\label{HP}
Suppose that $X$ is a subordinate Brownian motion satisfying {\bf (H)}.
For every $R>0$, there exists $c=c(R, \phi, d)>0$
such that for every $r\in (0
, R)$, every $x\in \R^d$, and every nonnegative function $h$ on
$\R^d$ which is harmonic in $B(x, r)$ with respect to $X$ we have
$$
\sup_{y \in B(x, r/2)}h(y) \le c \inf _{y \in B(x, r/2)}h(y).
$$
\end{thm}

Recall that $H=\{x=(\wt{x},x_d):\, \wt{x}\in \R^{d-1}, x_d>0\}$ is
the half-space in $\R^d$ and that $\partial H=\{x=(\wt{x},0):\,
\wt{x}\in \R^{d-1}\}$ denotes its boundary. Let $z\in \partial H$.
We will say that a function $h:\R^d\to \R$ vanishes continuously on
$ H^c \cap B(z, r)$ if $h=0$ on $ H^c \cap B(z, r)$ and $h$ is
continuous at every point of $\partial H\cap B(z,r)$.
By \cite[Lemma 4.2]{CKSV1} and its proof (which works for all
subordinate Brownian motions satisfying {\bf (H)}), we
 see that, if  $h$ is a nonnegative function in $\R^d$ that is
harmonic in $H\cap B(z, r)$ with respect to $X$ and vanishes
continuously on $H^c\cap B(z, r)$, then $h$ is regular harmonic in
$H\cap B(z, r)$ with respect to $X$.
Thus, using the Harnack inequality and a Harnack chain argument, the
following form of the boundary Harnack principle is a consequence of
\cite[Theorem 1.3]{KSV2}
and \cite[Theorem 4.22]{KSV3}
 in case $\alpha \in (0,2)$,
and \cite[Theorem 1.2]{KSV4} in case $\alpha=2$.
 Note that the distance of the point $x\in H$ to the boundary
$\partial H$ will be denoted by $\delta_H(x)$. Clearly,
$\delta_H(x)=x_d$ for $x \in H$.

\begin{thm}\label{t:bhp}
Suppose that $X$ is a subordinate Brownian motion satisfying {\bf (H)}.
Then, for every $R>0$ there exists a constant
$C_3(R)=C_3(d,\phi, R)>0$ such that for $r \in (0, 2R]$,
$z\in \partial H$ and any nonnegative function $h$ in $\R^d$ that is
harmonic in $H \cap B(z,r)$ with respect to $X$ and vanishes
continuously on $H^c \cap B(z, r)$, we have
\begin{equation}\label{e:bhp}
\frac{h(x)}{V(\delta_H(x))} \le C_3(R) \frac{h(y)}{V(\delta_H(y))}
\qquad \textrm{for every } x, y\in  H \cap B(z, r/2).
\end{equation}
\end{thm}

Recall that  $\tau_H=\inf\{t>0:\, X_t\notin H\}$.
The process $X^H=(X^H_t)_{t\ge 0}$ obtained by killing $X$ upon
exiting $H$ is defined by
$$
X^H_t:=\left\{\begin{array}{ll}
X_t\, ,& t<\tau_H\, ,\\
\partial\, ,& t\ge \tau_H\, ,
\end{array}\right.
$$
where $\partial$ is the cemetery point. The killed process $X^H$ is
a symmetric Hunt process.
Any function $h$ on $H$ is automatically extended to $\partial$ by
setting $h(\partial)=0$.
A nonnegative function $h$ is harmonic with respect to $X^H$ if for
every open set $B$ such that $B\subset \overline{B}\subset H$ and
every $x\in B$ it holds that $h(x)=\E_x[h(X^H_{\tau_B})]$. A
nonnegative function $s$ on $H$ is said to be excessive with respect
to $X^H$ if $s(x)\ge \E_x[s(X^H_{t})]$ for all $t>0$ and $x\in H$,
and $\lim_{t\downarrow 0}\E_x[s(X^H_{t})]=s(x)$ for all $x\in H$.
Since a subordinate Brownian motion is always a strong
Feller process, any function which is excessive with respect to
$X^H$ is lower semi-continuous.

For every $D \subset H$, $X$ admits a Green function $G^D:D\times
D\to (0,\infty]$ defined by
$$
G^D(x,y)=G(x,y)-\E_x[G(X_{\tau_D},y)]\, ,\qquad x,y\in D\, .
$$
The Green function $G^D$ is symmetric and continuous (in extended
sense). We extend the domain of the Green function $G^D$ to $\R^d
\times \R^d$ by setting $G^D$ to be zero outside of $D \times D$.

We will frequently use the well-known fact that $G^D(\cdot, y)$ is
harmonic in $D\setminus\{y\}$, and regular harmonic in $D\setminus
\overline{B(y,\varepsilon)}$ for every $\varepsilon >0$. Moreover,
by the strong Markov property, for all open sets $D_1 \subset D_2
\subset H$,
\begin{equation}\label{e:gsmp}
G^{D_1}(x,y)=G^{D_2}(x,y)-\E_x[G^{D_2}(X_{\tau_{D_1}},y)] \le
G^{D_2}(x,y),\quad x,y \in D_1.
\end{equation}

The following sharp estimates of $G^H$ for points close to the
boundary and to each other will be crucial in the sequel.

\begin{thm}\label{t:green}
Suppose that $X$ is a subordinate Brownian motion satisfying {\bf (H)}.
Then, for every $R>0$ there exists a constant
$C_4(R)=C_4(d,\phi, R)>0$  such that for all $x,y\in H$
satisfying $|x-y|<R$ and $\delta_H(x)\wedge \delta_H(y)<R$ it holds
that
\begin{eqnarray}
\lefteqn{C_4(R)^{-1} \left(1\wedge \frac{V(\delta_H(x))}{V(|x-y|)}\right)
\left(1\wedge \frac{V(\delta_H(y))}{V(|x-y|)}\right)G(x,y)\le G^H(x,y) }\nonumber \\
&\le & C_4(R)  \left(1\wedge \frac{V(\delta_H(x))}{V(|x-y|)}\right)
\left(1\wedge \frac{V(\delta_H(y))}{V(|x-y|)}\right)G(x,y)\, .\label{e:green}
\end{eqnarray}
\end{thm}

\pf Fix $R>0$. We first note that by the argument in
\cite[Lemma 5.1]{CKS2} (using \eqref{e:asymp-V-I} and \eqref{e:asymp-V-II}),
\eqref{e:green} is equivalent to
\begin{equation}\label{e:green2}
c^{-1}\left(1 \wedge\frac{V(\delta_H(x))V(\delta_H(y))}{V(|x-y|)^2} \right)
G(x,y) \le G^H(x,y)
\le c \left(1 \wedge\frac{V(\delta_H(x))V(\delta_H(y))}{V(|x-y|)^2} \right)
G(x,y)
\end{equation}
for some constant $c=c(d,\phi, R)>1$. By
translation invariance, we can assume $x,y \in B(0, \sqrt{5}R/2)
\cap H$.

Recall that a domain in $\bR^d$ is a connected open set in $\bR^d$.
A domain $D$ in $\bR^d$ is said to be a  $C^{1,1}$ domain
if there are $ R_0>0 $ and
$\Lambda_0>0$ such that for every $z\in\partial D$, there exist a
$C^{1,1}$-function $\varphi=\varphi_z: \bR^{d-1}\to \bR$ satisfying
$\varphi (0)=0, \nabla\varphi (0)=(0, \dots, 0)$, $| \nabla \varphi
(x)-\nabla \varphi (w)| \leq \Lambda_0 |x-w|$,
 and an orthonormal
coordinate system $CS_z$ : $y=(y_1, \dots, y_{d-1}, y_d):=(\wt y, \,
y_d)$  with origin at $z$ such  that $ B(z, R_0 )\cap D= \{y=(\wt y,
\, y_d)\in B(0, R_0) \mbox{ in } CS_z: y_d
> \varphi (\wt y) \}$.
The pair $(R_0, \Lambda_0)$ is called the $C^{1, 1}$ characteristics
of the $C^{1, 1}$ domain $D$.
Choose a \emph{bounded} $C^{1,1}$ domain $D$ such that $B(0, 3R)
\cap H \subset D \subset H$ and such that its $C^{1, 1}$
characteristics $(R_0, \Lambda_0)$ depends only on $d$ and $R$.
 The estimates \eqref{e:green2} with $H$ replaced by $D$ and for all
$x,y \in D$ were proved in \cite[Theorem 1.1]{KSV2} in case $\alpha\in (0,2)$,
and in \cite[Theorem 1.4]{KSV4} in case $\alpha=2$.
Note that $c$ depends on $d,\phi, R$ only.

It follows from \eqref{e:gsmp},
\cite[Theorem 1.1]{KSV2}  and \cite[Theorem 1.4]{KSV4}
that
$$
G^H(x,y) \ge G^D(x,y) \ge  c_1 \left(1 \wedge\frac{V(\delta_D(x))
V(\delta_D(y))}{V(|x-y|)^2} \right)
G(x,y)= c_2 \left(1 \wedge\frac{V(\delta_H(x))V(\delta_H(y))}
{V(|x-y|)^2} \right)
G(x,y),
$$
thus it suffices to show the second inequality in \eqref{e:green2}.

Let $x_R:=(\tilde 0, R/2)$. By Theorem \ref{t:bhp} applied to $G^H(\cdot, w)$,
\begin{align*}
&\int_{H\setminus D} G^H(y, w) \P_x(X_{\tau_D} \in dw)
\le C_3(R) \frac{V(\delta_H(y))}{V(R/2)} \int_{H\setminus D}
G^H(x_R, w)\P_x(X_{\tau_D} \in dw)\\
&\le c_3 V(\delta_H(y)) \int_{H\setminus D}  G(R)\P_x(X_{\tau_D} \in dw)
\le c_4 V(\delta_H(y)) \P_x(X_{\tau_D} \in H \setminus D).
\end{align*}
Using the boundary Harnack principle for $\P_{\cdot}(X_{\tau_D} \in
H \setminus D)$, we also have
\begin{eqnarray*}
\P_x(X_{\tau_D} \in H \setminus D) \le c_{5}
\P_{x_R}(X_{\tau_D} \in H \setminus D) \frac{V(\delta_H(x))}{V(R/2)}
\le c_{6} V(\delta_H(x)).
\end{eqnarray*}
Since
$$
V(\delta_H(x))V(\delta_H(y)) \le c_7 \left(1 \wedge\frac{V(\delta_H(x))
V(\delta_H(y))}{V(|x-y|)^2} \right) \le c_8  \left(1 \wedge
\frac{V(\delta_H(x))V(\delta_H(y))}{V(|x-y|)^2} \right)
G(x,y),
$$
we obtain
\begin{equation} \label{e:hgdfgh}
\int_{H\setminus D} G^H(y, w) \P_x(X_{\tau_D} \in dw)
\le c_9 \left(1 \wedge\frac{V(\delta_H(x))V(\delta_H(y))}{V(|x-y|)^2} \right)
G(x,y).
\end{equation}
Since,  by \eqref{e:gsmp} and the symmetry of $G^H$,
$$
G^H(x,y)=G^D(x,y)+\E_x[G^H(X_{\tau_D},y)]=G^D(x,y)+\E_x[G^H(y, X_{\tau_D})],
$$
 we obtain the second inequality in \eqref{e:green2}  using \eqref{e:hgdfgh},
\cite[Theorem 1.1]{KSV2} and \cite[Theorem 1.4]{KSV4}. \qed

\begin{remark}\label{r:green}{\rm
In case $\phi(\lambda)=a\lambda + b^\beta\lambda^{\beta/2}$ with $a, b>0$ and
$\beta\in (0, 2)$,
the estimates \eqref{e:green} are valid
 for all $x,y\in H$. This was proved in \cite[Theorem 1.7]{CKS2} by
using the heat kernel estimates obtained in \cite{CKS2} and a lengthy computation. }
\end{remark}

\section{Martin kernel and minimal thinness}

In this section we always assume that $X$ is a subordinate Brownian
motion satisfying ({\bf H}). Recall that $X^H$ is the process
obtained from $X$ by killing it upon exiting $H$.

In the remainder of this paper we will use $x_0$ to denote the point
$(\tilde{0}, 1)\in H$ and set
$$
M^H(x,y):=\frac{G^H(x,y)}{G^H(x_0, y)}\, ,\qquad x,y\in H, y\neq x_0\, .
$$
As the process $X^H$ satisfies Hypothesis (B) in \cite{KW}, $H$  has
a Martin boundary $\partial^MH$ with respect to
$X$ satisfying the
following properties:
\begin{description}
\item{(M1)} $H\cup \partial^M H$ is
compact metric space;
\item{(M2)} $H$ is open and dense in $H\cup \partial^M H$,  and its relative
topology coincides with its original topology;
\item{(M3)}  $M^H(x ,\, \cdot\,)$ can be uniquely extended  to $\partial^M H$
in such a way that,  $ M^H(x, y) $ converges to $M^H(x, w)$ as $y\to
w \in \partial^M H$, the function $x \to M^H(x, w)$  is excessive
with respect to $X^H$, the function $(x,w) \to M^H(x, w)$ is jointly
continuous on $H\times \partial^M H$ and $M^H(\cdot,w_1)\not=
M^H(\cdot, w_2)$ if $w_1 \not= w_2$;
\end{description}

For any $w\in \partial^M H$, the function $x \to M^H(x, w)$ is
called the Martin kernel of $X^H$ corresponding to $w$.
A point $w\in \partial^M H$ is called a finite Martin
boundary point if there is a bounded sequence $\{w_n\}\subset H$
converging to $w$ in the Martin topology.

Recall that a positive harmonic function $h$ with respect to $X^H$
is minimal if every positive harmonic function $g$ with respect to
$X^H$ such that $g\le h$ is proportional to
 $h$. The minimal Martin boundary of $ X^H$ is defined as
$$
\partial^mH=\{z\in \partial^MH: M^H(\cdot, z) \mbox{ is
minimal harmonic with respect to $X^H$ }\}.
$$
A point $z\in \partial^mH$ is called a finite minimal Martin
boundary point if there is a bounded sequence $\{w_n\}\subset H$
converging to $z$ in the Martin topology.

We will show that the finite part of the minimal Martin boundary of
$H$ with respect to
$X$ coincides with $\partial H$. The claim
above can be proved by using Theorems \ref{t:bhp}--\ref{t:green} and
following the methodology from \cite{B}, \cite[Section 5]{KSV1} and
\cite[Section 6]{CKSV2}. Since the Martin boundary of $H$ with
respect to $X^H$ contains also non-finite boundary points, slight
modifications of the argument are called for. In order to make the
paper more readable, we provide below the details for most of the
proofs.

The L\'evy measure of the process $X$ has a density $J$, called
the L\'evy density, given by
$$
J(x)=\int^{\infty}_0(4\pi t)^{-d/2}e^{-|x|^2/(4t)}\eta(t)dt, \qquad x\in \R^d.
$$
Recall that $\eta$ denotes the L\'evy density of the subordinator $S$.
Thus $J(x)=j(|x|)$ with
$$
j(r)
:=\int^{\infty}_0(4\pi t)^{-d/2}e^{-r^2/(4t)}\eta(t)dt, \qquad r>0.
$$
Note that the function $r\mapsto j(r)$ is continuous and decreasing on $(0, \infty)$.
When hypothesis {\bf (H)} is satisfied, $j$ enjoys the
following properties:
\begin{description}
\item{(a)} For any $K>0$, there exists $c=c(K)>0$ such that
\begin{equation}\label{H:1}
j(r)\le c\, j(2r), \qquad \forall r\in (0, K);
\end{equation}
\item{(b)} There exists $c>0$ such that
\begin{equation}\label{H:2}
j(r)\le c\, j(r+1), \qquad \forall r>1.
\end{equation}
\end{description}
(See \cite[(2.17)--(2.18)]{KSV2} in case $\alpha\in (0, 2)$ and \cite[(2.5)--(2.6)]{KSV4} in case $\alpha=2$.)
We further have (see \cite[Theorem 3.4]{KSV3}) that in case $\alpha\in (0,2)$,
\begin{equation}\label{e:estimate-of-j}
j(r)\asymp \frac{\phi(r^{-2})}{r^{d}}\, ,\quad r\to 0\, .
\end{equation}

For an open set $U\subset \R^d$, let
$$
K^U(x,z)\,:=\,
\int_U
{G^U(x,y)}
J(y-z)
 dy, \quad (x,z) \in U \times
\overline{U}^c.
$$
Then for any nonnegative measurable function $f$ on $\bR^d$,
$$
\E_x\left[f(X_{\tau_U});\,X_{\tau_U-} \not= X_{\tau_U}
\right] =\int_{\overline{U}^c} K^U(x,z)f(z)dz.
$$

Note that when subordinate Brownian motion $X$ satisfies {\bf (H)}
there exist  $c_1, c_2, r_1>0$ such
that for every $r \in (0, r_1]$ and $x_1 \in \R^d$,
\begin{eqnarray}
 K^{B(x_1,r)}(x,y) \,&\le &\, c_1 \, j(|y-x_1|-r) \left(\phi(r^{-2})\phi((r-|x-x_1|)^{-2})\right)^{-1/2} \label{P1}\\
 &\le &\, c_1 \, j(|y-x_1|-r) \phi(r^{-2})^{-1}\ \label{P1-worse}
\end{eqnarray}
for all $(x,y) \in B(x_1, r)\times \overline{B(x_1, r)}^c$ and
\begin{equation}\label{P2}
K^{B(x_1, r)}(x_1, y) \,\ge\, c_2\, j(|y-x_1|)
\phi(r^{-2})^{-1}, \qquad \textrm{ for all } y \in \overline{B(x_1, r)}^c.
\end{equation}
In case $\alpha\in (0,2)$, \eqref{P1} and \eqref{P2} were shown in \cite[Proposition 4.10]{KSV3}, while \eqref{P1-worse} follows from \eqref{P1} and the fact that $\phi$ is increasing. In case $\alpha=2$, the proof is analogous to the proof of \cite[Proposition 4.10]{KSV3}, except that one uses \cite[Proposition 3.1]{KSV09} instead of \cite[Proposition 4.9]{KSV3},  and \cite[Proposition 3.5]{RSV} instead of
\cite[Lemma 4.2]{KSV3}.

We use the notation that $A_r(Q):=(\tilde Q, r/2)$ for  $Q \in \partial H$.
Note that $H \cap B(Q,r)$ contains the ball $B(A_r(Q), r/2)$.

Using ({\bf H}) and \cite[Theorems
1.5.3 and 1.5.11]{BGT}, there exists a positive constant
$R_*:=R_*(d,\phi)<1 \wedge r_1$ such that
\begin{equation}\label{el8}
\int_{0}^r \frac{1}{s\phi(s^{-1})}ds \,\le \, 4 \,
\frac{1}{\phi(r^{-1})}, \qquad \forall\,  0<s<r\le 4 R_*,
\end{equation}
and for every $r \le 2 R_*$,
\begin{equation}\label{ll}
\frac1{2}\,\le\, \min \left(\frac{ \ell(8^{-2}
r^{-2})}{\ell(r^{-2})},\,
\frac{ \ell(16 r^{-2})} {\ell(r^{-2})}
\right)\,\le \, \max\left(\frac{\ell(8^{-2} r^{-2})}{\ell(r^{-2})},\,
\frac{ \ell(16 r^{-2})} {l(r^{-2})}
\right) \,\le\, 2.
\end{equation}
We will fix the constant $R_*$ in remainder of this section.

The next lemma is proved by the same argument as that of
\cite[Lemma 5]{B} and \cite[Lemma 5.2]{KSV1}
 in case $\alpha\in (0,2)$.

\begin{lemma}\label{l:5B}
Suppose that $X$ is a subordinate Brownian motion satisfying {\bf (H)}.
There exist positive constants  $c=c( d,\phi)$ and
$\gamma=\gamma(d,\phi)\in (0,\alpha)$
such that for all $Q\in \partial H$, $r\in (0, R_*]$, and nonnegative function $h$ in $\R^d$
which is harmonic with respect to $X$ in $H \cap B(Q, r)$ and vanishes continuously on $H^c\cap B(Q,r)$ we have
$$
\phi(r^{-2})\, h(A_r(Q))\le c 4^{-\gamma k}\phi(4^{2k} r^{-2})\, h(A_{4^{-k}r}(Q))\, , \qquad k=0, 1, \dots .
$$
In case $\alpha=2$ we may take $\gamma=1$.
\end{lemma}
\pf Without loss of generality, we may assume $Q=0$.
By Theorem \ref{t:bhp}, there exists a constant
$C_3(R_*)=C_3(d,\phi, R_*)>0$ such that for $r \in (0, 2R_*]$,
$$
h(A_r(0))\,\le\, C_3(R_*)\, \frac{V(\delta_H(A_r(0)))}
{V(\delta_H(A_{4^{-k}r}(0)))} h(A_{4^{-k}r}(0))
\,=\, C_3(R_*)\, \frac{V(r/2)}{V(4^{-k}r/2)} h(A_{4^{-k}r}(0)).
$$
This and \eqref{e:asymp-V-II} complete the proof  of the lemma in case $\alpha=2$.

Assume now that $\alpha\in (0,2)$. Fix $r \le  R_*$ and let $
\eta_k:=4^{-k}r$, $A_k:=A_{ \eta_k}(0)$ and $B_k\,:=\,B(A_k,
\eta_{k+1})$, $k=0,1, \dots. $ Note that the $B_k$'s are disjoint.
So by the harmonicity of $h$, we have
$$
h(A_k)
\,\ge\, \sum_{l=0}^{k-1} \E_{A_k}\left[h(X_{\tau_{B_k}}):\,
X_{\tau_{B_k}} \in B_l \right]\\
\,=\, \sum_{l=0}^{k-1} \int_{B_l} K^{B_k}(A_k, z) h(z) dz.
$$
Theorem  \ref{HP} implies that
$$
 \int_{B_l} K^{B_k}(A_k, z) h(z) dz \,\ge\, c_0\, h(A_l)
\int_{B_l} K^{B_k}(A_k, z) dz
$$
for some constant $c_0=c_0(d, \phi)>0$. Since dist$(A_k, B_l) \le
2 \eta_{l}$, by \eqref{P2}  and the monotonicity of $j$ we have
\begin{equation}\label{e:type1-estimate}
K^{B_k}(A_k,z)\ge\, c_1\, j(|2(A_k-z)|)\, \phi(\eta_{k+1}^{-2})^{-1}
\ge\, c_1\,  j(4\eta_l ) \, \phi(\eta_{k+1}^{-2})^{-1}\, ,
 \qquad
z \in B_l.
\end{equation}
Further, by \eqref{e:type1-estimate}, \eqref{H:1} and \eqref{e:estimate-of-j} we have
$$
K^{B_k}(A_k,z)\ge c_1 \frac{j(4\eta_l)}{\phi(\eta_{k+1}^{-2})}\ge c_2
\frac{j(\eta_{l+1})}{\phi(\eta_{k+1}^{-2})} \ge c_3
\frac{\phi(\eta_{l+1}^{-2})}{\eta_{l+1}^d \phi(\eta_{k+1}^{-2})}\, ,\quad z\in B_l
$$
for some constat $c_3=c_3(d,\phi)>0$.
Thus we have
$$
\int_{B_l}K^{B_k}(A_k, z)dz \ge\, c_4 \frac{\phi(\eta_{l+1}^{-2})}{\phi(\eta_{k+1}^{-2})}, \qquad z \in B_l
$$
for some constant $c_4=c_4(d, \phi)>0$. Therefore,
$$
h(A_k)\phi(\eta_{k+1}^{-2})\ge c_5 \sum_{l=0}^{k-1}h(A_l)\phi(\eta_{k+l}^{-2})
$$
for some constant $c_5=c_5(d, \phi)>0$.
Let $a_k :=  h(A_k)\phi(\eta_{k+1}^{-2})$ so that
$a_k \ge  c_5\sum_{l=0}^{k-1}  a_l$. By
induction, one can easily check that $ a_k  \ge c_6 (1+c_5/2)^{k} a_0$
for some constant $c_6=c_6(d, \alpha)>0$.
Thus, with $ \gamma ={\ln(1+\frac{c_5}2)} (\ln 4)^{-1}>0 $ we get
$$
\phi(r^{-2})h(A_r(Q))\,\le\, c_6\, 4^{-\gamma k}\,
\phi(4^{2k}r^{-2}) h(A_{4^{-k}r}(Q)).
$$
Note that we can choose $c_5>0$ small enough so that $\gamma <\alpha$.
This completes the proof in case $\alpha\in (0,2)$.
\qed

By modifying the proof of \cite[Lemma 5.3]{KSV1},
one can obtain the following
\begin{lemma}\label{l:la}
Suppose that $X$ is a subordinate Brownian motion satisfying {\bf (H)}.
Suppose $Q \in \partial H$ and $r \in (0, R_*]$. If $w
\in H\setminus B(Q, r)$, then
$$
G^H(A_r(Q), w)\, \ge\, c\, \phi(r^{-2})^{-1} \int_{B(Q, r)^c} J (z- Q) G^H(z, w)dz
$$
for some constant $c=c( d,\phi)>0$.
\end{lemma}

\pf
This proof is similar that of \cite[Lemma 5.3]{KSV1}, we give the details here for completeness.
Without loss of generality, we may assume $Q=0$. Fix $w \in
H\setminus B(0, r)$ and let $A:=A_r(0)$ and $h(\cdot) := G^H(\cdot, w)$.
Since $h$ is regular harmonic in $H\cap B(0,
3r/4)$ with respect to $X$, we have
\begin{eqnarray*}
h(A) &\ge& \E_A \left[ h\left( X_{\tau_{H \cap
B(0,3r/4)}}\right); X_{\tau_{H \cap B(0,3r/4)}}
\in B(0, r)^c \right]= \int_{B(0, r)^c}K^{H\cap
B(0, 3r/4)}(A, z)h(z)dz\\
&=&\int_{B(0, r)^c} \int_{H\cap B(0,  3r/4 ) }
 G^{H\cap B(0,3r/4 )}(A,y)\,J(y-z) dyh(z)dz.
\end{eqnarray*}
Since $B(A, r/4)\subset H \cap B(0, 3r/4 )$, by the
monotonicity of the Green functions,
$$
G^{H\cap B(0,3r/4 )}(A,y) \, \ge \,
G^{B(A, r/4)}(A,y),
\quad y \in B(A, r/4).
$$
Thus
$$
h(A) \ge  \int_{B(0, r)^c }
\int_{ B(A, r/4)  }
G^{B(A, r/4)}(A,y) J(y-z) dyh(z)dz = \int_{B(0, r)^c }
K^{B(A, r/4)}(A, z)h(z)dz,
$$
which is greater than or equal to $ c_2\phi((r/4)^{-2})^{-1}\int_{B(0, r)^c}
J(z-A)h(z)\,dz$ for some positive constant $c_2=c_2(d, \phi)$ by \eqref{P2}.
Now the conclusion of the lemma follows immediately from \eqref{ll} and
\eqref{H:1}--\eqref{H:2}.
\qed

Using the above lemma,
\eqref{H:1}, \eqref{H:2} and \eqref{P1-worse},
the
proof of next lemma is  almost the same as that of \cite[Lemma 5.4]{KSV1}.

\begin{lemma}\label{l:14B}
Suppose that $X$ is a subordinate Brownian motion satisfying {\bf (H)}.
There exist positive constants
$c_1=c_1(d,\phi)$ and $c_2=c_2(d,\phi)<1$ such that for any
$Q \in \partial H $, $r\in (0, R_*]$ and  $w \in H \setminus
B(Q,4r)$,  we have
$$
\E_x\left[G^H(X_{\tau_{H \cap B_k}}, w):\,X_{\tau_{H
\cap B_k}} \in B(Q, r)^c
 \right] \,\le\, c_1\,c_2^{k} \, G^H(x,w), \quad x \in H \cap B_k,
$$
where $B_k:=B(Q, 4^{-k}r)$, $ k=0,1, \dots$.
\end{lemma}

\pf Without loss of generality, we may assume $Q=0$. Fix $r \le R_*$
and  $w \in H \setminus B(0,4r)$.  Let $\eta_k:=4^{-k}r $,
$B_k:=B(0,\eta_k)$ and
$$
u_k(x)  \,:=\, \E_x\left[G^H(X_{\tau_{H \cap B_k}}, w);
X_{\tau_{H \cap B_k}}
 \in B(0, r)^c \right], \quad x \in H \cap B_k.
$$
Note that for $x\in H\cap B_{k+1}$
\begin{eqnarray}\label{e:dec1}
u_{k+1}(x)
=   \E_x\left[ G^H(X_{\tau_{H \cap B_k}}, w)  ;\, \tau_{H
\cap B_{k+1}} = \tau_{H\cap B_{k}}, ~ X_{\tau_{H \cap B_{k}}}
\in B(0, r)^c  \right]  \le  u_{k}(x).
\end{eqnarray}
Let $A_k\,:=\,A_{\eta_k}(0) $. Since $G^H(\,\cdot\,, w)$ is zero on
$H^c$, we have
\begin{eqnarray*}
&&u_{k}(A_k) \,= \, \E_{A_k}\left[G^H(X_{\tau_{H \cap B_k}}, w);\,
X_{\tau_{H \cap B_{k}}}
 \in B(0, r)^c  \right] \\
&&\le  \E_{A_k}\left[G^H(X_{\tau_{ B_k}}, w)  ;\, X_{\tau_{
B_{k}}}
 \in B(0, r)^c \right]
\,=\,  \int_{B(0, r)^c} K^{B_k}(A_k,z) G^H(z,w) dz.
\end{eqnarray*}
By \eqref{P1-worse} we have that there exists a constant
$c_1=c_1(d, \phi)>0$ such that
\begin{equation}\label{e:dec2}
u_{k}(A_k) \,\le\, c_1\,
\phi(\eta_k^{-2})^{-1} \int_{B(0, r)^c}j(|z|-\eta_k) G^H(z,w)\, dz
\quad k=1,2, \dots.
\end{equation}
 Note that
 $
|z|/2 \le |z|-\eta_k \le |z|
 $ for $z \in
B(0, 4R_*) \setminus B(0,r)$ and $
|z|-R_* \le |z|-\eta_k \le |z|
 $  for $z \in
B(0, 4R_*)^c$. Thus  using \eqref{H:1}--\eqref{H:2}, we see that
Lemma \ref{l:la} and (\ref{e:dec2})  imply that
$$
u_{k}(A_k) \,\le\, c_2
\frac{\phi(r^{-2})}{\phi((4^{-k}r)^{-2})}\, G^H(A_0,w)\, ,
\quad k=1,2, \dots.
$$
Now applying Lemma \ref{l:5B}  we get
$$
u_{k}(A_k) \,\le\, c_3
4^{-\gamma k} G^H(A_k,w)\, ,
\quad k=1,2, \dots.
$$
Finally,
Theorem \ref{t:bhp} gives that for $k=1,2, \dots$
$$
\frac{u_k(x)}{G^H(x, w)}    \,\le\,
\frac{u_{k-1}(x)}{G^H(x, w)}
\,\le\, c_4\, \frac{u_{k-1}(A_{k-1})}{G^H(A_{k-1}, w)}  \,\le\,
c_5\,
4^{-\gamma k}
.
$$
\qed

Now the next theorem follows from Theorem \ref{t:bhp} and Lemma
\ref{l:14B} in the same way as \cite[Lemma 16]{B}
follows from \cite[Lemmas 13--14]
{B}. We omit the details.

\begin{thm}\label{t:martin}
Suppose that $X$ is a subordinate Brownian motion satisfying {\bf (H)}.
For each $x\in H$ and each $z\in \partial H$ there
exists the limit
$$
M^H(x,z):=\lim_{y\to z}M^H(x,y)\, .
$$
Moreover, the mapping
$(x, z)\mapsto M^H(x, z)$ is continuous on
$H\times\partial H$.
\end{thm}

\begin{remark}\label{r:martin}
Note that Theorem \ref{t:martin} shows that the finite part of the
Martin boundary of $H$ can be identified with a subset of $\partial
H$.
\end{remark}

We recall here that the Martin kernel for the killed $\alpha$-stable
process in $H$ is explicitly known (see \cite{Bog}) and is given by
\begin{eqnarray*}
M^H(x,z)&=&\frac{\delta_H(x)^{\alpha/2}}{|x-z|^d}\, (1+|z|^2)^{d/2}\, ,
\quad x\in H,\,   z\in \partial H\, ,\\
M^H(x,\infty)&=&\delta_H(x)^{\alpha/2}\, ,\quad x\in H\, .
\end{eqnarray*}
The same form of the Martin kernel is valid for the killed Brownian
motion in $H$ with $\alpha$ replaced by 2. We cannot hope to obtain
explicit formulae for the Martin kernel for the process $X^H$, but
the following sharp two-sided estimates for
points close to the origin and to each other
will suffice for our
purpose.

\begin{thm}\label{t:martin-estimate}
Suppose that $X$ is a subordinate Brownian motion satisfying {\bf (H)}.
For every $R>0$, there exists a constant
$C_5(R)=C_5(d,\phi, R)>0$
 such that for all $x\in H$ and all
$z\in \partial H$ satisfying $|z|<R$ and $|x-z|< R/2$ it holds that
$$
C_5(R)^{-1}\frac{V(\delta_H(x))}{|x-z|^d}\, (1+|z|^2)^{d/2}\le M^H(x,z)
\le C_5(R) \frac{V(\delta_H(x))}{|x-z|^d}\, (1+|z|^2)^{d/2}\, .
$$
\end{thm}
\pf Let $x\in H$ and $z\in \partial H$ such that $|z|<R$ and
$|x-z|<R/2$. Let $y\in H$ be such that $|y-z|<R/2$ and $\delta_H(y)
\le \delta_H(x) \wedge |x-y|$. Then $|x-y|<R$, $\delta_H(x)\wedge
\delta_H(y )<R$ and $\delta_H(x)\vee \delta_H(y)<2|x-y|$. Hence, by
Theorem \ref{t:green} and \eqref{e:green2}
$$
\frac{G^H(x,y)}{G^H(x_0,y)}\asymp \frac{\frac{V(\delta_H(x))
V(\delta_H(y))}{V(|x-y|)^2}\, G(x,y)}{\frac{V(\delta_H(x_0))
V(\delta_H(y))}{V(|x_0-y|)^2}\, G(x_0,y)}=\frac{V(\delta_H(x))}
{V(1)} \frac{G(x,y)}{V(|x-y|)^2}\frac{V(|x_0-y|)^2}{G(x_0,y)}.
$$
Let $y\to z$; by Theorem \ref{t:martin} the left-hand side converges
to $M^H(x,z)$, while the right-hand side converges to
$$
\frac{V(\delta_H(x))}{V(1)} \frac{G(x,z)}{V(|x-z|)^2}
\frac{V(|x_0-z|)^2}{G(x_0,z)}\, .
$$
By use of Proposition \ref{p:GoverV}(i) and the fact that
$|x_0-z|=(1+|z|^2)^{1/2}$ it follows that
$$
M^H(x,z)\asymp \frac{V(\delta_H(x))}{|x-z|^d}\, (1+|z|^2)^{d/2}\, ,
$$
which proves the theorem.
\qed

Theorem \ref{t:martin-estimate} in particular implies that
$M^H(\cdot , z_1)$ differs from $M^H(\cdot, z_2)$ if $z_1$ and $z_2$
are two different  points on $\partial H$. Together with Remark
\ref{r:martin}, this shows that the finite part of the Martin
boundary of $H$ can be identified with $\partial H$.

Now using Theorem \ref{t:green}, \eqref{P1} and Lemma
\ref{l:5B}, we can prove the next two lemmas by the arguments in the
proofs of \cite[Lemmas 4.6--4.7]{KS2} and \cite[Lemmas
5.6--5.7]{KSV1}.

\begin{lemma}\label{l:MH1}
Suppose that $X$ is a subordinate Brownian motion satisfying {\bf (H)}.
For every $z \in \partial H$ and $B \subset \overline{B}
\subset H$, $M^H(X_{\tau_{B}} , z)$ is $\P_x$-integrable.
\end{lemma}

\pf Take a sequence $\{z_m\}_{m \ge 1} \subset H\setminus
\overline{B}$ converging to $z$. Since $M^H (\cdot , z_m )$ is
regular harmonic with respect to  $X$ in $B$, by Fatou's lemma and
Theorem \ref{t:martin},
$$  \E_x \left[ M^H\left(X_{\tau_{B}} , z\right)\right]\,=\,
\E_x \left[ \lim_{m\to \infty} M^H\left(X_{\tau_{B}} ,
z_m\right) \right]\,\leq\, \liminf_{m\to \infty} M^H(x, z_m)\,=\,
M^H(x,z) \,<\,\infty .
$$
\qed

\begin{lemma}\label{l:MH2}
Suppose that $X$ is a subordinate Brownian motion satisfying {\bf (H)}.
For every $z \in \partial H$ and  $x \in H$,
$$
M^H(x,z) \,=\, \E_x \left[M^H\left(X^H_{\tau_{B(x,r)}} ,
z\right)\right], \quad \mbox{ for
 every }
  0<r \le R_* \wedge \frac13\delta_H(x).
$$
\end{lemma}

\pf Without loss of generality we assume that $z=0$, and fix $x \in
H$ and $r \le R_* \wedge \frac13\delta_H(x)$. Let $
\eta_m\,:=4^{-m}r$ and  $z_m\,:=\, A_{ \eta_m}(0),$  $m=0,1, \dots$.
By the harmonicity of $M^H(\cdot , z_m )$, to prove the lemma, it  suffices to show that
$\{ M^H(X_{\tau_{B(x,r)}}, z_m ) : m\geq m_0 \}$ is $\P_x $-uniformly
integrable for some large $m_0$.

Using Theorem \ref{t:bhp}, there exist constants $m_0
\ge 0$ and $c_1>0$ such that for every $w \in H \setminus B (z,
\eta_m)$,
\begin{equation}\label{eqn:4.2}
M^H(w, z_m) \,=\, \frac{G^H(w, z_m)}{G^H(x_0 , z_m )} \,\leq\,
c_1\, \lim_{y \to z} \frac{G^H(w, y)}{G^H(x_0 , y )} \,=\, c_1\,M^H (w, z),
\quad m\geq m_0.
\end{equation}
Using (\ref{eqn:4.2}) and Lemma \ref{l:MH1}, for
any $\eps>0$, there is an $N_0>1$ such that
\begin{eqnarray}\label{eqn:c4.3}
\lefteqn{ \E_x \left[M^H\left(X_{\tau_{B(x,r)}}, z_m\right) ; \,
M^H\left(X_{\tau_{B(x,r)}}, z_m\right) > N_0 ~\mbox{ and }~
X_{\tau_{B(x,r)}} \in H \setminus B(z, \eta_m)\right]} \nonumber \\
&\leq & c_1\, \E_x \left[M^H\left(X_{\tau_{B(x,r)}}, z\right) ;
\, c_1 M^H\left(X_{\tau_{B(x,r)}}, z\right) > N_0 \right]
 \,<\,   c_1\,\frac{ \eps}{4c_1} \,=\, \frac{\eps}{4} .
\end{eqnarray}
On the other hand, by \eqref{P1-worse}, we have for $m\geq m_0$,
\begin{eqnarray*}
&&\E_x \left[M^H\left(X^H_{\tau_{B(x,r)}}, z_m\right) ; \,
X_{\tau_{B(x,r)}} \in H
 \cap B(z, \eta_m)\right]\,=\,\int_{H \cap B(z, \eta_m)} M^H(w, z_m)
K^{B(x,r)}(x,w) dw \\
&&\leq \, c_2
\int_{H \cap B(z, \eta_m)} M^H(w, z_m) j(|w-x|-r) \phi(r^{-2})^{-1} \, dw
\end{eqnarray*}
for some $c_2=c_2(d, \phi)>0$.
Since $|w-x| \ge |x-z|-|z-w|\ge \delta_H(x)- \eta_m \ge 3r -r=2r$,
using the monotonicity of $j$ in the above equation, we see that
\begin{eqnarray}
&&\E_x \left[M^H\left(X^H_{\tau_{B(x,r)}}, z_m\right) ; \,
X_{\tau_{B(x,r)}} \in H
 \cap B(z, \eta_m)\right]\nonumber\\
&&\leq \, c_3 j(r)\phi(r^{-2})^{-1}\, \int_{B(z, \eta_m)} M^H(w, z_m) dw
\le c_4 \,G^H(x_0, z_m)^{-1} \int_{B(z,  \eta_m)} G^H(w, z_m)
dw\label{e:new1}
\end{eqnarray}
for some  $c_3=c_3(\phi)>0$ and $c_4=
c_4(\phi, r)>0$.
By Lemma \ref{l:5B} there exists there exist
$c_5=c_5(\phi, m_0, r)$ and $\gamma\in (0,\alpha)$ such that
\begin{equation}\label{e:auxilliary-1}
G^H(x_0,z_m)^{-1}\le c_5 4^{-\gamma(m-m_0)} \frac{\phi(4^{2m}r^{-2})}{\phi(4^{2 m_0} r^{-2})}\, G^H(x_0,z_{m_0})^{-1}\, .
\end{equation}
Further, by \eqref{e:asymp-G-I}  we have that there are constants
$c_i=c_i(\phi, m_0, r)>0$, $i=6,7,8,9$, such that
\begin{eqnarray}
\int_{B(0,\eta_m)} G^H(w,z_m)\, dw &\le & c_6 \int_{B(z_m,2\eta_m)}\frac{dw}{|w-z_m|^d \phi(|w-z_m|^{-2})} \nonumber \\
&\le & c_7 \int_0^{2\eta_m} \frac{ds}{s\phi(s^{-1})} \le c_8 \phi((2\eta_m)^{-2})^{-1}
\,\le\,  c_9 \phi(4^{2m}r^{-2})^{-1}\, ,\label{e:auxilliary-2}
\end{eqnarray}
where the third inequality follows from \eqref{el8} and the fourth from \eqref{ll}.
Putting together \eqref{e:auxilliary-1} and \eqref{e:auxilliary-2} we arrive at
\begin{equation}\label{e:auxilliary-3}
G^H(x_0,z_m)^{-1}\int_{B(z,\eta_m)} G^H(w,z_m) \le c_{10}4^{-\gamma(m-m_0)} \frac{G^H(x_0,z_{m_0})}{\phi(4^{2 m_0}r^{-2})}
\end{equation}
for a constant $c_{10}=c_{10}(\phi, m_0, r)>0$.
Using \eqref{eqn:c4.3},
\eqref{e:new1} and \eqref{e:auxilliary-3},
we can take $m_1=m_1(\eps, \phi, m_0, r) \ge m_0$ large enough so that  for $m\geq m_1$,
\begin{eqnarray*}
\lefteqn{\E_x  \left[M^H\left(X_{\tau_{B(x,r)}}, z_m\right);\,
M^H\left(X_{\tau_{B(x,r)}}, z_m\right) > N_0\right]}\\
&\leq&  \E_x  \left[M^H\left(X_{\tau_{B(x,r)}}, z_m\right) ; \,
X_{\tau_{B(x,r)}} \in H \cap B(z, \eta_m)\right]\\
&&+ \E_x  \left[M^H\left(X_{\tau_{B(x,r)}}, z_m\right) ; \,
M^H\left(X_{\tau_{B(x,r)}}, z_m\right) > N_0 \,\mbox{ and }\,
X_{\tau_{B(x,r)}} \in H \setminus B(z, \eta_m)\right] \\
&\le & c_4 c_{10} 4^{-\gamma(m-m_0)} \frac{G^H(x_0,z_{m_0})}{\phi(4^{2 m_0}r^{-2})} +\frac{\epsilon}{4}
\end{eqnarray*}
is less than $\epsilon$. As each $M^H(X_{\tau_{B(x,r)}}, z_m)$ is $\P_x$-integrable, we conclude that $\{ M^H(X_{\tau_{B(x,r)}}, z_m ) : m\geq m_0 \}$  is uniformly integrable under $\P_x$.
\qed

The proof of next result is taken from \cite[Theorem 6.10]{CKSV2}.

\begin{thm}\label{T:L4.3}
Suppose that $X$ is a subordinate Brownian motion satisfying {\bf (H)}.
For every $z \in \partial H$,
the  function $x\mapsto M^H(\cdot, z)$ is harmonic in $H$
with respect to
$X^H$.
\end{thm}
\pf
Fix $z\in \partial H$ and let $h(x):=M^H(x,z)$. Consider an
open set $D_1 \subset\overline{D_1}\subset H$, $x \in D_1$, and put
$r(x)=
R_* \wedge \frac13\delta_H(x)$ and $ B(x)=B(x,r(x))$.
Define a sequence of stopping times $\{T_m, m\ge1\}$ as follows:
$
T_1:=\inf\{t>0: X_t\notin B(X_0)\},
$
and for $m\ge 2$,
$$
T_m:= \begin{cases}  T_{m-1}+
\tau_{B(X_{T_{m-1}})}\circ\theta_{T_{m-1}}
\qquad &\hbox{if }X_{T_{m-1}}\in D_1\\
\tau_{D_1} \qquad &\hbox{otherwise}.
\end{cases}
$$
Note  that $X^H_{\tau_{D_1}}\in \partial D_1$ on
$\cap_{n=1}^{\infty}\{T_n < \tau_{D_1}\}. $ Thus, since $\lim_{m
\to  \infty} T_m=\tau_{D_1}$ $\P_x$-a.s. and $h$ is continuous in
$D$,
using the quasi-left-continuity of $X^H$,
we have
$
\lim_{m \to  \infty} h(
X^H_{T_m})= h(
X^H_{\tau_{D_1}})$ on $\cap_{n=1}^{\infty}\{T_n < \tau_{D_1}\}.
$
Now, by the dominated
convergence theorem and  Lemma \ref{l:MH2},
\begin{align*}
&h(x)
= \lim_{m\to\infty}\E_x\left[h(
X^H_{T_m});\, \cup_{n=0}^{\infty}
\{T_n=\tau_{D_1} \}\right] + \lim_{m\to\infty} \E_x
\left[h(
X^H_{T_m});\, \cap_{n=0}^{\infty}
\{T_n<\tau_{D_1} \}\right]\\
&= \E_x\left[h(
X^H_{\tau_{D_1}});\,
\cup_{n=0}^{\infty}
\{T_n=\tau_{D_1} \}\right] +  \E_x\left[h(
X^H_{\tau_{D_1}});\,
\cap_{n=0}^{\infty}
\{T_n<\tau_{D_1} \}\right] =  \E_x\left[h(
X^H_{\tau_{D_1}})\right].
\end{align*}
\qed

The following two lemmas will serve as a counterpart of Theorem
\ref{t:martin-estimate} for estimating the Martin kernel $M^H(x,w)$
when $x$ and $w$ are far away.

\begin{lemma}\label{l:boy12}
Suppose that $X$ is a subordinate Brownian motion satisfying {\bf (H)}.
For every $R>0$ and $\eps \in (0,1)$, there exists a
constant $c(R, \eps)=c(d,\phi, R, \epsilon)>0$ such
that for all $z\in \partial H$ with $|z|<R$, all $x\in H\cap
B(z,\eps)$, and all $w\in \partial^M H\setminus (\partial H\cap
B(z,2 \eps))$ it holds that
$$
c^{-1}V(\delta_H(x))\le M^H(x,w)\le c V(\delta_H(x))\, .
$$
\end{lemma}
\pf Fix $z\in\partial H$ with $|z|<R$ and $w \in \partial^M H
\setminus (\partial H\cap B(z,2\eps))$. Let $y_0:=(\tilde z,
\eps/2)$. Consider a sequence $\{ w_n\} \subset H \cap B(z, \eps/2)$
such that $w_n \to w$. Then, by Theorem \ref{t:bhp},
$$
C_3(\eps)^{-1} \frac{V(\delta_H(x))}{V(\eps/2)}\le
\frac{G^H(x, w_n)}{G^H(y_0, w_n)} \, \le \, C_3(\eps)
\frac{V(\delta_H(x))}{V(\eps/2)}, \quad \text{ for every }
n \ge 1
\text{ and } x \in H \cap B(z, \eps).
$$
Now, using Theorem \ref{HP}, we see that
$$
c^{-1}\,C_3(\eps)^{-1} \frac{V(\delta_H(x))}{V(\eps/2)}\le
\frac{G^H(x, w_n)}{G^H(x_0, w_n)}  \,\le\, c\,C_3(\eps)
\frac{V(\delta_H(x))}{V(\eps/2)}, \quad \text{ for every }
n \ge 1 \text{ and } x \in H \cap B(z, \eps)
$$
where $c=c(R, \phi, d)>0$. Thus, letting $n \to \infty$,
we conclude that
$$
c^{-1} \, C_3(R)^{-1} \frac{V(\delta_H(x))}{V(1/2)} \le
M^H(x, w)\,\le \, c \, C_3(R) \frac{V(\delta_H(x))}{V(1/2)}
\quad \text{ for every }  x \in H \cap B(z, \eps),
$$
which finishes the proof.
\qed

\begin{lemma}\label{l:boy13}
Suppose that $X$ is a subordinate Brownian motion satisfying {\bf (H)}.
For every $R>0$, there exists a constant
$c=c(d,\phi, R)>0$ such that for all $z\in
 \partial H$ with $|z|<R$ and all $x\in H$ such that $|x-z|>R/2$, it
holds that
$
M^H(x,z)\,\le \,c\, G(|x-z|+1).
$
\end{lemma}

\pf Fix $z\in\partial H$ with $|z|<R$, $x\in H$ such that
$|x-z|>R/2$, and let $z_0:=(\tilde z, (R \wedge 1)/4)$. Consider a
sequence $(z_n)_{n\ge 1}\subset H \cap B(z, (R \wedge 1)/4)$ such
that $z_n\to z$. By Theorem \ref{t:bhp} applied to function
$G^H(x,\cdot)$ and $G^H(x_0,\cdot)$ we have that
$$
M^H(x,z) = \lim_{n \to \infty} \frac{G^H(x,z_n)}{G^H(x_0,z_n)}\le C_3(
(R \wedge 1)/4)) \frac{G^H(x,z_0)}{G^H(x_0,z_0)}\, .
$$
It follows from Theorem \ref{t:green}, the fact that $|x_0-z_0|<R+1$
and the monotonicity of $G$ and $V$ that
$$
G^H(x_0,z_0)\ge C_4(R+1)^{-1}G(R+1)\left(1\wedge
\frac{V(1)}{V(R+1)}\right)\left(1\wedge \frac{V( (R \wedge
1)/4))}{V(R+1)}\right) \ge c_1\, ,
$$
where $c_1=c_1(d,\phi, R)>0$. Hence, by using that
$G^H(x,z_0)\le G(x,z_0)=G(|x-z_0|) \le G(|x-z|+1)$, we
prove the lemma.
\qed

\begin{remark}\label{r:martin-goes-to-zero}{\rm
By combining Theorem \ref{t:martin-estimate} and Lemma \ref{l:boy13}
we see that for every $z\in \partial H$ and every $w\in
\partial^M H\setminus \{z\}$, it holds that $\lim_{x\to
w}M^H(x,z)=0$. Indeed, if $w\in \partial H$, then there exists $R>1$
large enough such that $|z|<R$ and $|z-w|<R/3$ and the claim follows
from Theorem \ref{t:martin-estimate}. On the other hand, if $w\in
\partial^M H\setminus \partial H$, then we use Lemma \ref{l:boy13}
and the fact that $\lim_{|x|\to \infty}G(|x-z|+1)=0$. }
\end{remark}

With the explicit estimates from Theorem \ref{t:martin-estimate} and
Lemmas \ref{l:boy12} and \ref{l:boy13}, a slight modification of the
argument in \cite[Theorem 3.7]{CS} gives the following.

\begin{thm}\label{t3.1}
Suppose that $X$ is a subordinate Brownian motion satisfying {\bf (H)}.
For every $z\in \partial H$, the function $M^H(\cdot ,z)$
is  minimal harmonic in $H$ with respect to
$X$.
Therefore,
 the finite part of the minimal Martin boundary of $H$ of $X^H$
coincides with $\partial H$.
\end{thm}
\pf Fix $z\in \partial H$ and suppose that $h$ is nonnegative
function harmonic with respect to $X^H$ satisfying $h\le M^H(\cdot,
z)$. By \cite{KW}, there is a measure $\mu$ on $\partial^M H$ such
that
$$
h(x)=\int_{\partial^M H} M^{H}(x,w)\, \mu(dw)\, .
$$
We want to prove that $\mu$ is a multiple of the point mass at $z$.
Let $\delta\in (0,1)$ and denote by $\nu$ the restriction of $\mu$
to $\partial^M H \setminus (\partial H\cap B(z,\delta))$. If we show
that $\nu= 0$, we are done (since $\delta>0$ is arbitrary). Let
\begin{equation}\label{e:representation-u}
u(x):= \int_{\partial^M H}M^H(x,w)\, \nu(dw)\, .
\end{equation}
Then $u$ is a nonnegative harmonic functions with respect to $X^H$
such that $u(x)\le h(x)\le M^H (x,z)$ for all $x\in H$. Let
$\epsilon \in (0,1)$. It follows from Theorem
\ref{t:martin-estimate} and Lemma \ref{l:boy13} (as explained in
Remark \ref{r:martin-goes-to-zero}) that for every $w\in \partial^M
H \setminus (\partial H\cap B(z,\epsilon))$ we have that $\lim_{x\to
w} M^H(x,z)=0$. Hence, $\lim_{x\to w} u(x)=0$. Suppose,
additionally, that $ 4\epsilon <\delta$.
It follows from Proposition \ref{l:boy12} that
$$
\lim_{x\to \partial H\cap B(z, 2\epsilon)} M^H(x,w)=0\qquad
\textrm{uniformly in }w\in \partial^M H\setminus (\partial H\cap
B(z,\delta))\, .
$$
Applying the dominated convergence theorem in
 \eqref{e:representation-u}, we conclude that $\lim_{x\to
w} u(x)=0$ for every $w\in \partial H\cap B(z,2\epsilon)$ as well.
This shows that $u$ is a nonnegative harmonic function with respect
to $X^H$ that vanishes continuously at $\partial^M H$. Therefore it
is identically equal to zero proving that its representation measure
$\nu=0$. \qed

For simplicity we use the notation $\partial_\infty^m H:=\partial^m
H \setminus \partial H$. It is known (e.g., \cite{KW}) that every
function $s$ which is excessive with respect to $X^H$
can be uniquely represented in the form
\begin{equation}\label{e:representation}
s(x)=G^H\mu (x)+M^H\nu(x):=\int_H G^H(x,y)\, \mu(dy) +
\int_{\partial^m H}M^H(x,z)\, \nu(dz)\, , \quad x\in H\, ,
\end{equation}
where $\mu$ is a measure on $H$ and $\nu$ a finite measure on
$\partial^m H$.

For $z\in \partial H$, let $X^{H,z}=(X^H_t, \P_x^z)$ be the
$M^H(\cdot, z)$-process.
The existence of such a process is discussed in \cite{Fol, KW}. Let
$\zeta$ denote the lifetime of $X^{H,z}$.
It is known (see \cite{KW}) that $\lim_{t\uparrow
\zeta-}X^{H,z}_t=z$
a.s. $\P_x^z$.

For $A\subset H$, let $T_A:=\inf\{t> 0:\, X^H_t\in A\}$. A Borel
set $A\subset H$ is said to be \emph{minimally thin} in $H$ with
respect to $X$ at $z\in \partial H$ if there exists $x\in H$ such
that $\P^z_x(T_A <\zeta)\neq 1$, i.e., if with positive probability
the conditioned process $X^{H,z}$ starting  from $x$ does not hit
$A$. It is proved in \cite[Satz 2.6]{Fol}  that $A$ is minimally
thin in $H$ with respect to
$X$
 at $z$ if and only if
$$
P_A M^H(\cdot, z)\neq M^H(\cdot, z)\, .
$$
Here $P_A$ denotes the hitting operator to $A$ for $X^H$: $P_A
f(x)=\E_x[f(X^H_{T_A})]$. In
potential-theoretic language, $P_A M^H(\cdot, z)=
\wh{R}^A_{M^H(\cdot, z)}$ -- the balayage of $M^H(\cdot, z)$ onto
$A$. By following the proof of \cite[Theorem 9.2.6]{AG} and
using \cite[Lemma 2.7]{Fol} instead of \cite[Lemma 9.2.2(c)]{AG}
 one can show that $A$ is minimally thin in $H$ with respect to
$X$ at
$z\in \partial H$ if and only if there exists a function $s$
excessive with respect to $X^H$ such that
\begin{equation}\label{e:thinness}
\liminf_{x\to z, x\in A}\frac{s(x)}{M^H(x,z)} > \nu(\{z\})\, ,
\end{equation}
where $\nu$ is the representing measure of the harmonic part $M^H
\nu$  of $s$ (as in \eqref{e:representation}).
By subtracting
$\nu(\{z\})$ from the measure $\nu$, we have

\begin{prop}\label{p:new}
Suppose $A\subset H$ and  $z\in \partial H$. The following are
equivalent:

\noindent (a)
$A$ is minimally thin in $H$ with respect to
$X$ at $z$;

\noindent (b) There
exists an excessive function $s(x)=G^H\mu (x)+M^H\nu(x)$
such that \eqref{e:thinness} holds;

\noindent (c) There
exists an excessive function $s'(x)=G^H\mu (x)+M^H\nu'(x)$
 with $\nu'(\{z\})=0$ such that
$$
\liminf_{x\to z, x\in A}\frac{s'(x)}{M^H(x,z)}>0\, .
$$
\end{prop}

\section{Main result}

The key to the proof of our generalization of
Theorems \ref{t:dahlberg}--\ref{t:burdzy}
is the following lemma and proposition, which are generalizations of
\cite[Theorem 2]{Sj2}.

\begin{lemma}\label{l:sjogren}
Suppose that $X$ is a subordinate Brownian motion satisfying {\bf (H)}.
Let $\nu$ be a finite measure on $\partial^m H$ with
$\nu(\{0\})=0$. Then
$$
\int_{\{x \in H: M^{H}\nu (x) >M^H(x, 0)/2\}\cap B(0,1)}
|x|^{-d}\, dx<\infty\, .
$$
\end{lemma}
\pf
Note that since $\nu(\{0\})=0$,
$\nu(\{z\in \partial H: |z|<\epsilon\})$
can be
made arbitrarily small by choosing $\epsilon $ small enough. Let
$\epsilon <1/4$ be
a positive constant such that
\begin{equation}\label{e00}
\nu(\{z\in \partial H:\, |z|<\epsilon\})<\frac{1}{10}\, C_5(1)^{-2}
 \left(\frac{{\sqrt 3}}{8}\right)^{d} 2^{-d/2}.
\end{equation}

Set
$A:=\{x\in H: M^H\nu(x) >M^H(x, 0)/2\}\cap B(0,1)$
and  let
$$
\nu^{(1)}=\nu|_{\{z\in \partial H:\, |z| < \epsilon\}}\, ,
\quad \nu^{(2)}=\nu|_{\{z\in \partial H:\, \epsilon\le |z|<1\}}\, ,
\quad \nu^{(3)}=\nu|_{\{z\in \partial H:\, |z|\ge 1\}\cup \partial ^m_\infty H}
$$
so that $\nu=\nu^{(1)}+\nu^{(2)}+\nu^{(3)}$.

Let $y_0=(\tilde{0}, 1/8)$ and $x\in H$ be a point with $|x|<\kappa$
where $\kappa<1/4$ is to be chosen later. For $z \in {\{z\in
\partial H:\, |z|\ge 1\}\cup \partial ^m_\infty H}$ the function
$M^H(\cdot, z)$ is harmonic and vanishes continuously on $H^c\cap
B(0,1/2)$ by Theorem \ref{t:martin-estimate} and
Remark \ref{r:martin-goes-to-zero}.
 Hence by Theorem \ref{t:bhp} it holds that
$$
\frac{M^H(x,z)}{M^H(y_0,z)}\le C_3(1/4)
\frac{V(\delta_H(x))}{V(\delta_H(y_0))}\, .
$$
Therefore, by Theorem \ref{t:martin-estimate},
 for $x\in H$ with $|x|<\kappa$,
\begin{eqnarray*}
M^H \nu^{(3)}(x)
&\le & C_3(1/4) \frac{V(\delta_H(x))}{V(1/8)}
\int_{\{z\in
\partial H:\, |z|\ge 1\}\cup \partial ^m_\infty H}M^H(y_0,z)\, \nu(dz)\\
&\le& C_3(1/4)\frac{\kappa^d V(\delta_H(x))}{V(1/8)} |x|^{-d}
\int_{\{z\in \partial H:\, |z|\ge 1\}\cup \partial ^m_\infty H}
M^H(y_0,z)\, \nu(dz) \\
&\le&
 C_3(1/4)\frac{\kappa^dC_5(1)}{V(1/8)}M^H(x,0)
\int_{\{z\in \partial H:\, |z|\ge 1\}\cup \partial ^m_\infty H}
M^H(y_0,z)\, \nu(dz).
\end{eqnarray*}
Now choose $\kappa< 1/4$ small enough so that
$$
 C_3(1/4)\frac{\kappa^dC_5(1)}{V(1/8)}\int_{\{z\in \partial H:\, |z|\ge
1\}\cup \partial ^m_\infty H} M^H(y_0,z)\, \nu(dz)<1/10.
$$
Then for $x\in H$ with $|x|<\kappa$,
\begin{equation}\label{e0}
M^H \nu^{(3)}(x)< \frac{1}{10} M^H(x,0)\, .
\end{equation}

Let $0<\rho<1$, $x\in H$ with $|x|<\rho \epsilon$ and $z\in \partial
H$ with $\epsilon \le |z|<1$. Then we have $2>|x-z|\ge
|z|-|x|>(1-\rho)\epsilon$. Thus, by Theorem \ref{t:martin-estimate},
\begin{eqnarray*}
M^H \nu^{(2)}(x)&=&\int_{\{z\in \partial H: \epsilon\le |z|<1 \}}M^H(x,z)\, \nu(dz)\\
&\le &C_5(4)\, V(\delta_H(x))\int_{\{z\in \partial H: \epsilon\le
|z|<1\}}
|x-z|^{-d}(1+|z|^2)^{d/2}\, \nu(dz)\\
&\le &C_5(4)\, 2^{d/2}V(\delta_H(x)) (1-\rho)^{-d}
\epsilon^{-d}\nu(\{z\in \partial H: \epsilon\le |z|<1\})\\
&\le &C_5(4)\, 2^{d/2}
\nu(\partial^m H) V(\delta_H(x))(1-\rho)^{-d}
\rho^d |x|^{-d}\\
&\le &C_5(4)^2\, 2^{d/2} \nu(\partial^m H) \left(\frac{\rho}{1-\rho}
\right)^d M^H(x,0)\, .
\end{eqnarray*}
Now choose $\rho$ small enough so that $C_5(4)^2\, 2^{d/2}
\nu(\partial^m H) \rho^d(1-\rho)^{-d}<1/10$. Then for $x\in H$ with
$|x|<\rho \epsilon$,
\begin{equation}\label{e1}
M^H \nu^{(2)}(x)< \frac{1}{10} M^H(x,0)\, .
\end{equation}

Consider now $M^H\nu^{(1)}(x)$ for $x\in H$ with $|x| < \rho \eps
\wedge \kappa$. If $\delta_H(x)>|x|/2$, then $|x|\le 2\delta_H(x)\le
2|x-z|$ for each $z\in \partial H$. Thus  by Theorem
\ref{t:martin-estimate},
\begin{align*}
M^H(x,z) \, \le \, &C_5(1)\, V(\delta_H(x)) |x-z|^{-d}(1+|z|^2)^{d/2}\\
\, \le \,  & C_5(1)\, V(\delta_H(x)) (|x|/2)^{-d}2^{d/2}
\, \le \,  C_5(1)^2 \, 2^{3d/2} M^H(x,0)\, ,
\end{align*}
for all $z\in \partial H$ with $|z|<\epsilon$. Hence, for $x\in H$
with $\delta_H(x)>|x|/2$ we have
\begin{eqnarray}\label{e2}
M^H\nu^{(1)}(x)&=&\int_{\{z\in \partial H: |z|<\epsilon\}}M^H(x,z)\,
\nu^{(1)}(dz)\nonumber\\
&\le& C_5(1)^2 2^{3d/2} M^H(x,0)\nu(\{z\in
\partial H: |z|<\epsilon\}) \, <\, \frac{1}{10}M^H(x,0)
\end{eqnarray}
by the choice of $\epsilon$ in \eqref{e00}.

It remains to study $M^H\nu^{(1)}$ on the set $E:=\{x\in H:\,
\delta_H(x)<|x|/2\}\cap B(0, \rho\eps \wedge \kappa)$. Note that for
$x \in E$ we have $|x| < \frac{2}{\sqrt 3} |\tilde{x}|$.
For $j\ge 0$, set
$$
R_j:=\{z\in \partial H:\, 2^{-j-1}<|z|\le 2^{-j}\}\quad \text{ and
}\quad \wt{R}_j:=\{x\in E:\, \wt{x}\in R_j\}.
$$
Put $\nu_j=\nu^{(1)}|_{R_j}$.
Assume that $x\in \wt{R}_j$ and $z\in R_k$ for $|j-k|>1$. If $k
\ge 2+j$, $|x-z| \ge |\tilde{x}|-|z| \ge 2^{-j-1}-2^{-k} \ge \frac14
2^{-j}\ge \frac14 |\tilde{x}| \ge \frac{{\sqrt 3}}{8}|x|$. If $j \ge
k+2$, $|x-z| \ge |z| -\frac{2}{\sqrt 3}  |\tilde{x}|\ge
2^{-k-1}-\frac{2}{\sqrt 3} 2^{-j} \ge 2^{-j+1}-\frac{2}{\sqrt 3}
2^{-j} \ge (2-\frac{2}{\sqrt 3}) |\tilde{x}| \ge \frac{{\sqrt
3}}{2}(2-\frac{2}{\sqrt 3})|x|=({\sqrt 3}-1)|x|$.
Thus, if $|j-k|>1$, we have $|x-z|>\frac{{\sqrt 3}}{8}|x|$ for every
$x\in \wt{R}_j$ and $z\in R_k$, and so  by Theorem
\ref{t:martin-estimate}, for $x\in \wt{R}_j$,
\begin{align*}
M^H\nu_k(x) \, = \, &\int_{R_k}M^H(x,z)\, \nu^{(1)}(dz)
\, \le \, C_5(1) \int_{R_k}V(\delta_H(x))|x-z|^{-d}(1+|z|^2)^{d/2}\, \nu^{(1)}(dz)\\
\, \le \, & C_5(1)\, (\frac{8}{{\sqrt 3}})^d 2^{d/2}
V(\delta_H(x))|x|^{-d}\nu^{(1)}(R_k) \, \le \, C_5(1)^2\,
(\frac{8}{{\sqrt 3}})^d 2^{d/2} \nu^{(1)}(R_k) M^H(x,0)\, ,
\end{align*}
implying
\begin{equation}\label{e3}
\sum_{|k-j|>1}M^H \nu_k(x)\le C_5(1)^2\, (\frac{8}{{\sqrt 3}})^d 2^{d/2}
\nu(\{|z|<\epsilon\})  M^H(x,0) \le \frac{1}{10}M^H(x,0)
\end{equation}
by the choice of $\epsilon$ in \eqref{e00}.

Let
$$
M_j:=\{x\in \wt{R}_j:\, M^H(\nu_{j-1}+\nu_j+\nu_{j+1})(x)>
\frac25 M^H(x,0)\}.
$$
Note that by Theorem \ref{t:martin-estimate},
\begin{eqnarray*}
M_j &\subset &\{x\in \wt{R}_j:\, \int |x-z|^{-d}(1+|z|^2)^{d/2}
(\nu_{j-1}(dz)+\nu_j(dz)+\nu_{j+1}(dz))>\frac{2}{5C_5(1)^2}
|x|^{-d}\}\\
&\subset &\{x\in \wt{R}_j:\, \int |x-z|^{-d}2^{d/2}(\nu_{j-1}(dz)+
\nu_j(dz)+\nu_{j+1}(dz))>\frac{2}{5C_5(1)^2}|x|^{-d}\}\\
&=&\{x\in \wt{R}_j:\, |x|^{-d}\ast (\nu_{j-1}+\nu_j+\nu_{j+1}) >
\frac{2^{1-d/2}}{5C_5(1)^2}|x|^{-d}\}\, .
\end{eqnarray*}
For $x\in \wt{R}_j$ we have that $2^{-j-1}<|x|< \frac{2}{\sqrt 3}
2^{-j}$. Therefore,
$$
M_j\subset \{|x|^{-d}\ast (\nu_{j-1}+\nu_j+\nu_{j+1}) >
c_1^{-1} 2^{jd} \}
$$
where $c_1:=2^{3d/2-1} 3^{-d/2}5 C_5(1)^{2}$.
It follows from \cite[Lemma 1 and Theorem]{Sj1} that
$$
|M_j| \, \le \,  c_1 2^{-jd}(\nu^{(1)}
(R_{j-1})+\nu^{(1)}(R_{j})+\nu^{(1)}(R_{j+1})).
$$
Therefore by using again that $2^{-j-1}<|x|< \frac{2}{\sqrt 3}
2^{-j}$, we get that
$$
\int_{M_j}|x|^{-d}\, dx \, \le \,  c_2 2^{jd}|M_j|\le c_3
(\nu^{(1)}(R_{j-1})+\nu^{(1)}(R_{j})+\nu^{(1)}(R_{j+1}))
$$
implying
\begin{equation}\label{e4}
\sum_j \int_{M_j}|x|^{-d}\, dx \, \le \,  c_3 \sum_j
(\nu^{(1)}(R_{j-1})
+\nu^{(1)}(R_{j})+\nu^{(1)}(R_{j+1}))
 \, \le \,  3 c_3 \nu(\partial^m H)\, .
\end{equation}
It follows from \eqref{e0}--\eqref{e4} that
\begin{eqnarray*}
&&\int_{\{x \in H, |x|<\rho\eps \wedge \kappa  :    M^H\nu (x)  >
M^H(x,0)
/2 \}}|x|^{-d} \, dx\\
&\le & \int_{\{x \in H, |x|<\rho\eps \wedge \kappa  :
M^H\nu^{(1)} (x)  > 3 M^H(x,0)/10 \}}|x|^{-d} \, dx \\
&= & \int_{\{x \in E:    M^H\nu^{(1)} (x)  > 3 M^H(x,0)/10 \}}
|x|^{-d} \, dx \\
&= & \sum_j  \int_{\{x \in  \tilde R_j :  \sum_{|j-k|>1}
M^H\nu_k (x)  +\sum_{|j-k| \le 1}   M^H\nu_k (x)  > 3
M^H(x,0)/10 \}}|x|^{-d} \, dx \\
&\le& \sum_j \int_{M_j}|x|^{-d}\, dx\le 3 c_3
\nu(\partial^m H) \, <\,  \infty.
\end{eqnarray*}
Therefore, using the trivial fact that $\int_{\{\rho\eps \wedge
\kappa <|x|<1\}} |x|^{-d}\, dx <\infty$, we conclude that
$$
\int_{A}|x|^{-d} \, dx \, \le \, \sum_j \int_{M_j}|x|^{-d}\, dx +
\int_{\{ \rho\eps \wedge \kappa <|x|<1\}}
|x|^{-d}\, dx \, <\, \infty\, .
$$
\qed

\begin{prop}\label{p:sjogren}
Suppose that $X$ is a subordinate Brownian motion satisfying {\bf (H)}.
Let $s=G^H\mu+M^H\nu$ be an excessive function
with respect to $X^H$
represented as \eqref{e:representation}
with $\nu(\{0\})=0$. Let $A:=\{x\in H:\,
s(x)>M^H(x,0)\}$. Then
$$
\int_{A\cap B(0,1)}|x|^{-d}\, dx<\infty\, .
$$
\end{prop}

\pf First note that $A\subset \{G^H\mu >M^H(\cdot, 0)/2\}\cup
\{M^H\nu >M^H(\cdot, 0)/2\}$. By Lemma \ref{l:sjogren}, we only need
to consider $A_1:=\{G^H\mu >M^H(0,\cdot)/2\}\cap B(0,1)$.

Let $y_0=(\tilde{0}, 1/2)$ and $\kappa \ge 3$. By the boundary
Harnack principle (Theorem \ref{t:bhp}) applied to $G^H(\cdot, y)$
we obtain that,
$$
\int_{H\cap B(0,\kappa)^c}G^H(x,y)\, \mu(dy) \, \le \, C_3(1)
\frac{V(\delta_H(x))}{V(\delta_H(y_0))}\int_{H\cap
B(0,\kappa)^c} G^H(y_0,y)  \mu(dy)
$$
for every $\kappa \ge 3$ and $x \in B(0,1) \cap H$.
Now choose $\kappa$ large enough so that
\begin{equation}\label{fixkappa}
\frac{C_3(1)}{V(\delta_H(y_0))}\int_{H\cap B(0,\kappa)^c} G^H(y_0,y)
\mu(dy) < \frac1{4C_5(1)},
\end{equation}
where $C_5(1)$ is the constant from Theorem \ref{t:martin-estimate}.
Hence, for $x\in B(0, 1)\cap H$, we have
\begin{equation}\label{s3}
\int_{H\cap B(0,\kappa)^c}G^H(x,y)\, \mu(dy) \, \le \,
\frac1{4C_5(1)} V(\delta_H(x)) \, \le \, \frac1{4C_5(1)}
\frac{V(\delta_H(x))}{|x|^d} \, \le \,  \frac14 M^H(x,0)\, .
\end{equation}

Let ${\mathcal Q}=\{Q=(\prod_{i=1}^{d-1}  (k_i 2^{-j}, (k_i+1)2^{-j}
] ) \times (2^{-j}, 2^{-j+1} ]: k_i \in \bZ, j \in \bN\}$. Then
${\mathcal Q}$ is a cover of $H \cap \{x_d \le 1\}$. Clearly  the
cubes $Q$ are disjoint and $d^{-1/2} \textrm{diam}(Q) \le
\delta_H(Q) \le \textrm{diam}(Q).$ For $Q=(\prod_{i=1}^{d-1}  (k_i
2^{-j}, (k_i+1)2^{-j} ] ) \times (2^{-j}, 2^{-j+1} ]\in \mathcal Q$,
let $Q^*:=(\prod_{i=1}^{d-1}  ((k_i-1) 2^{-j}, (k_i+2)2^{-j} ]
)\times (2^{-j-1}, 2^{-j+2} ]$. Then ${\mathcal Q}^*=\{Q^*:\,  Q\in
{\mathcal Q}\}$ is a cover of $H \cap\{x_d \le 2\}$, $Q\subset
Q^*\subset H$ for each $Q\in {\mathcal Q}$, and
$2\delta_H(Q^*)=\delta_H(Q)$. Moreover, there exists a positive
integer $K=K(d)$ such that each cube in ${\mathcal Q}^*$ intersects
at most $K$ other cubes from  ${\mathcal Q}^*$. Let ${\mathcal Q}_1$
be the collection of cubes $Q\in \mathcal Q$ such that $H\cap B(0, 1)\cap Q$
is non-empty. For $x\in H\cap B(0,1)$,  let $Q(x)$
be the cube $Q$ in
${\mathcal Q}_1$ containing $x$ and denote by $Q^*(x)$ the
corresponding $Q^*$.

Let $G^H\mu=s_1+s_2+s_3$ where
\begin{eqnarray*}
s_1(x)&:=&\int_{Q^*(x)}G^H(x,y)\, \mu(dy)\, , \\
s_2(x)&:=&\int_{(H\setminus Q^*(x))\cap B(0,\kappa)}
G^H(x,y)\, \mu(dy)\, ,\\
s_3(x)&:=&\int_{(H\setminus Q^*(x))\cap B(0,\kappa)^c}
G^H(x,y)\, \mu(dy)\, ,
\end{eqnarray*}
where $\kappa \ge 3$ is the fixed constant in \eqref{fixkappa}.

Define
 $\gamma(r):=\int_{B(0,r)}G(0,z)\, dz$, $r>0$. Since $G^H(z,y)\le
G(z,y)$ we have for any $Q\in {\mathcal Q}_1$,
\begin{eqnarray*}
\int_{Q} s_1(w)\, dw &=& \int_{Q} \int_{Q^*} G^H(w,y)\, \mu(dy)\, dw
= \int_{Q^*}\mu(dy)\int_{Q} G^H(w,y)\, dw \\
&\le &\int_{Q^*}\mu(dy)\int_{Q} G(w,y)\, dw \le
\int_{Q^*}\mu(dy)\int_{B(y,3\mathrm{diam}(Q))}
G(w,y)\, dw\\
&\le & \mu(Q^*)\gamma(3\mathrm{diam}(Q))\le \mu(Q^*) \gamma(3 {\sqrt
d}\, \delta_H(Q))\, .
\end{eqnarray*}
Since $|w|\le (1+\sqrt{d})\textrm{dist}(0,Q)$ for $w\in Q$, we have
by Theorem \ref{t:martin-estimate},
$$
M^H(w,0) \, \ge \,  C_5(1)^{-1}V(\delta_H(w))|w|^{-d} \, \ge \, C_5(1)^{-1}
(1+\sqrt{d})^{-d}V(\delta_H(Q))\textrm{dist}(0,Q)^{-d} \quad \text{
for }w\in Q.
$$
We need to estimate the Lebesgue measure of
$$
B:=\{w\in Q:\, s_1(w)> C_5(1)^{-1} (1+\sqrt{d})^{-d}V(\delta_H(Q))\,
\textrm{dist}(0,Q)^{-d}\}.
$$
Since
$$
\int_Q s_1(w)\, dw \, \ge \,  C_5(1)^{-1}(1+\sqrt{d})^{-d} V(\delta_H(Q))\,
\textrm{dist}(0,Q)^{-d} |B|,
$$
we have that
\begin{eqnarray*}
|B|&\le &  C_5(1)(1+\sqrt{d})^dV(\delta_H(Q))^{-1}\,
\textrm{dist}(0,Q)^d
\int_Q s_1(w)\, dw \\
&\le & C_5(1)(1+\sqrt{d})^d  V(\delta_H(Q))^{-1}\,
\textrm{dist}(0,Q)^d \mu(Q^*)
\gamma(   3 {\sqrt d}\, \delta_H(Q))\\
&\le &C_5(1)(1+\sqrt{d})^d  C_4(3 {\sqrt d}) V(\delta_H(Q))^{-1} V(3
{\sqrt d} \delta_H(Q))^2\, \textrm{dist}(0,Q)^d\, \mu(Q^*),
\end{eqnarray*}
where in the third line we used Proposition \ref{p:GoverV}(ii).
Then
\begin{eqnarray*}
&&\int_{\{s_1>M^H(\cdot,0)/8\}\cap B(0, 1)}|w|^{-d}\, dw\\
&=& \sum_{Q \in {\mathcal Q}_1}\int_{\{s_1>M^H(\cdot,0)/8\}\cap
B(0, 1)\cap Q}|w|^{-d}\, dw\\
&\le & \sum_{Q \in {\mathcal Q}_1} \int_{\{s_1>M^H(\cdot,0)/8\}
\cap B(0, 1)\cap Q} \mathrm{dist}(0,Q)^{-d}\, dw\\
&\le&  \sum_{Q \in {\mathcal Q}_1} \int_{\{s_1>
C_5(1)^{-1}(1+\sqrt{d})^{-d}
V(\delta_H(Q))\mathrm{dist}(0,Q)^{-d}/8\}
\cap B(0, 1)\cap Q}\textrm{dist}(0,Q)^{-d}\, dw\\
&\le & C_5(1)  C_2(3 {\sqrt d})(1+\sqrt{d})^d  \sum_{Q \in {\mathcal Q}_1}
V(\delta_H(Q))^{-1} V(3 {\sqrt d}\delta_H(Q))^2\, \mu(Q^*)\\
&= & C_5(1)  C_2(3 {\sqrt d})(1+\sqrt{d})^d  \sum_{Q \in
{\mathcal Q}_1} V(2\delta_H(Q^*))^{-1} V(6 {\sqrt d}\delta_H(Q^*))^2\, \mu(Q^*)\\
&\le &c C_5(1)  C_2(2 {\sqrt d})(1+\sqrt{d})^d  \sum_{Q \in
{\mathcal Q}_1} V(\delta_H(Q^*)) \, \mu(Q^*)\, ,
\end{eqnarray*}
where the last line follows from \eqref{e:asymp-V-I} and
\eqref{e:asymp-V-II}. Note that
 $Q^*\subset B(0,4)$ for every $Q\in {\mathcal Q}_1$ and clearly
$\delta_H(Q^*) \le \delta_H(y)$ for all $y\in Q^*$. Thus by first
using these observations and the fact that each $Q^*\in {\mathcal
Q}^*$ intersects at most $K$ other cubes from ${\mathcal Q}^*$, and
then using Theorem \ref{t:bhp}, we conclude that the last sum is
dominated by
\begin{eqnarray*}
(K+1) \int_{H\cap B(0, 4)} V(\delta_H(y))\, \mu(dy)&\le
&(K+1)\frac{C_3(4)
V(2)}{G^H(w_1,y_1)} \int_{H\cap B(0,4)} G^H(w_1,y)\, \mu(dy)\\
&\le &(K+1) \frac{C_3(4)V(2)}{G^H(w_1,y_1)}G^H\mu(w_1)<\infty\, ,
\end{eqnarray*}
where $w_1=(\wt{0},10)$ and $y_1=(\wt{0},2)$.
 Therefore,
\begin{equation}\label{e6}
\int_{\{s_1>M^H (\cdot,0)/8\}\cap B(0, 1)}|w|^{-d}\, dw < \infty\, .
\end{equation}

Consider now $x \in \{s_2>M^H(\cdot, 0)/4\}$. For
$y\in H\setminus Q^*(x)$ it holds that
$|x-y|\ge 4^{-1} \delta_H(x)$,
and hence
$$
s_2(x)\le \int_{\{\delta_H(x)\le 4|x-y|\}
 \cap B(0,\kappa)} G^H(x,y)\, \mu(dy)\, .
$$
Let $\psi:H\to \partial H\setminus\{0\}$ be
a Borel function such that
$|y-\psi(y)|\le 2\delta_H(y)$ and define a measure $\mu'$ on
$\partial H\setminus\{0\}$ by
$$
\int_{\partial H\setminus\{0\}} f \, d\mu' =c_1
\int_{H\cap B(0,\kappa)} (f\circ \psi )(y)V(\delta_H(y))\, \mu(dy)\,
,
$$
where $c_1:=C_4(1+\kappa) C_1(1+\kappa)(11)^dC_5(11(1+\kappa))$.
Then by Theorem \ref{t:bhp},
$$
\mu'(\partial H\setminus\{0\})=c_1\int_{H\cap B(0,\kappa)}
V(\delta_H(y))\, \mu(dy)\le
c_1\frac{C_3(\kappa)
V(1)}{G^H(v_1, v_2)}G^H\mu(v_1)<\infty
$$
where $v_1=(\wt{0}, 3\kappa)$ and $v_2=(\wt{0},1)$. Hence, $\mu'$ is
positive and bounded.

We claim that
\begin{equation}\label{e:new}
s_2(x) \le M^H \mu'(x)  \quad \text{if } x \in \{s_2>M^H(\cdot,
0)/8\}\cap B(0, 1).
\end{equation}

Suppose that $y \in B(0, \kappa)\cap H$ and $x\in
\{s_2>M^H(\cdot,0)/8\}\cap B(0, 1)$ with
$\delta_H(x)\le 4|x-y|$.
Then
\begin{eqnarray*}
&&|x-\psi(y)|\le |x-y|+|y-\psi(y)|\le
|x-y|+2\delta_H(y)\\
&&\le |x-y|+2\delta_H(x)+2 |x-y|\le 11 |x-y| \le 11 (1+ \kappa) .
\end{eqnarray*}
Thus by Theorem \ref{t:martin-estimate},
\begin{eqnarray*}
M^H(x,\psi(y))&\ge& C_5(
11(1+ \kappa))^{-1}
V(\delta_H(x))(1+|\psi(y)|^2)^{d/2}|x-\psi(y)|^{-d}\\
&\ge& \frac{1}{(11)^d C_5(11 (1+ \kappa))} V(\delta_H(x))|x-y|^{-d}\, .
\end{eqnarray*}
Note that by Theorem \ref{t:green} and Proposition \ref{p:GoverV}(i),
$$
G^H(x,y) \, \le \,  C_4(1+\kappa) \frac{V(\delta_H(x))
V(\delta_H(y))}{V(|x-y|)^2}G(x,y) \, \le \, C_4(1+\kappa)
C_1(1+\kappa) V(\delta_H(x)) V(\delta_H(y))|x-y|^{-d}\, .
$$
Hence we get
$$
G^H(x,y) \, \le \,  c_{1} V(\delta_H(y)) M^H(x,\psi(y)).
$$
Therefore, for $x \in \{s_2>M^H(\cdot,
0)/8\}\cap B(0, 1),$
\begin{eqnarray*}
s_2(x)&\le & \int_{\{\delta_H(x)\le 2|x-y|\}
\cap B(0,\kappa)} G^H(x,y)\, \mu(dy)\\
&\le& c_{1} \int_{\{\delta_H(x)\le 2|x-y|\} \cap B(0,\kappa)}
V(\delta_H(y)) M^H(x,\psi(y))\, \mu(dy)\\
&\le &c_{1} \int_{H\cap B(0,\kappa)} V(\delta_H(y))
M^H(x,\psi(y))\, \mu(dy)\\
& =&\int M^H(x,z)\, \mu'(dz) \, = \, M^H\mu'(x)\, .
\end{eqnarray*}
Hence, we have proved \eqref {e:new}, and so, by Lemma \ref{l:sjogren},
\begin{equation}\label{e7}
\int_{\{s_2>M^H (\cdot,0)/8\}\cap B(0,1)}|x|^{-d}\, dx \, \le \,
\int_{\{M^H \mu' > M^H (\cdot,0)/8\}\cap B(0,1)} |x|^{-d}\, dx
\, <\, \infty.
\end{equation}

Together with \eqref{s3}, (\ref{e6}) and (\ref{e7}) we obtain that
$$
\int_{A_1}|x|^{-d}\, dx \, \le \,  \int_{\{s_1(x)+ s_2(x)>M^H(\cdot,0)/4\}
\cap B(0,1)}|x|^{-d} dx \, \le \, \sum_{i=1}^2\int_{\{s_i>M^H
(\cdot,0)/8\}\cap B(0,1)} |x|^{-d} dx \, < \, \infty.
$$
\qed

\begin{thm}\label{t:main1}
Suppose that $X$ is a subordinate Brownian motion satisfying {\bf (H)}.
Let $A$ be a Borel subset of $H$ and assume that
\eqref{c:criterion1} holds. Then $A$ is not minimally thin in $H$
with respect to
$X$ at $z=0$.
\end{thm}

\pf Assume that $A$ is minimally thin in $H$ with respect to
$X$
at 0. Then, using Proposition \ref{p:new} and multiplying $s$ by a
constant if necessary, there exists an excessive function
$s=G^H\mu+M^H \nu$ with $\nu(\{0\})=0$ such that
$$
\liminf_{x\to 0, x\in A}\frac{s(x)}{M^H(x,0)}=2\, .
$$

Let $B:=\{x\in H:\, s(x)>M^H (x,0)\}$.
Using the lower semi-continuity of $s$, we get that $B$ is an open
set. Thus there exists $\epsilon >0$
 such that $A\cap B(0,\epsilon)\subset B$. By Proposition
\ref{p:sjogren},
$$
\int_{A\cap B(0,\epsilon)}|x|^{-d}\, dx \,  \le \, \int_B|x|^{-d}\,
dx \, < \, \infty\, .
$$
Since $\int_{\{\epsilon <|x|<1\}}|x|^{-d}\, dx <\infty$, we have
proved the theorem. \qed

\begin{thm}\label{t:main2}
Suppose that $X$ is a subordinate Brownian motion satisfying {\bf (H)}.
Let $f:\R^{d-1}\to [0,\infty)$ be a Lipschitz function
with Lipschitz constant $a>0$. The set $A:=\{x=(\wt{x},x_d)\in H:\,
0<x_d\le f(\wt{x})\}$ is minimally thin in $H$ with respect to
$X$
at $z=0$ if and only if \eqref{c:criterion2} holds.
\end{thm}

\pf One direction follows by the same argument as the one in
\cite[Lemma 1]{Gar}. In fact, assume that the integral in
(\ref{c:criterion2}) diverges. For $x\in A=\{x=(\wt{x},x_d)\in H:\,
0<x_d\le f(\wt{x})\}$ we have
$$
|x| \, = \, (|\wt{x}|^2+x_d^2)^{1/2} \, \le \,
(|\wt{x}|^2+f(\wt{x})^2)^{1/2} \, \le \, (|\wt{x}|^2+a^2
|\wt{x}|^2)^{1/2} \, = \, |\wt{x}|(1+a^2)^{1/2}\, .
$$
Hence,
$$
\int_{A\cap B(0,1)} \frac{1}{|x|^d}\, dx\ge \int_{\{|\wt{x}|<1\}}
\int_0^{f(\wt{x})} \frac{1}{|\wt{x}|^d (1+a^2)^{1/2}}\,
dx=(1+a^2)^{-1/2} \int_{\{|\wt{x}|<1\}} f(\wt{x}) |\wt{x}|^{-d}\,
d\wt{x}= \infty\, .
$$
By Theorem \ref{t:dahlberg}, $A$ is not minimally thin in $H$ with
respect to
$X$ at 0.

The proof of the other direction is modeled after the proof of
\cite[Theorem 9.7.1]{AG}. Suppose that (\ref{c:criterion2}) holds
true.  Then
$$
V(\delta_H(x))\int_{\{|\wt{z}|<1\}} |x-\wt{z}|^{-d}|
\wt{z}|^{-d}f(\wt{z})\, d\wt{z} =V(\delta_H(x))\int_{\partial H\cap
B(0,1)} |x-z|^{-d} |z|^{-d} f(z)\, d\wt{\lambda}(z)<\infty
$$
where $\wt{\lambda}$ denotes the $d-1$-dimensional Lebesgue measure.
Define $s:H\to [0,\infty)$ by
\begin{equation}
s(x)=\int_{\partial H\cap B(0,1)}  M^H(x,z)(1+|z|^2)^{-d/2}
|z|^{-d} f(z)\, d\wt{\lambda}(z)\, .
\end{equation}
Then by Theorem \ref{t:martin-estimate},
$$
s(x) \, \le \,  C_5(4)\, V(\delta_H(x))\int_{\partial H\cap B(0,1)}
|x-z|^{-d} |z|^{-d} f(z)\, d\wt{\lambda}(z)\, \quad
\text{for every }x \in B(0,1) \cap H.
$$
Therefore, $s$ is well-defined (finite), and is harmonic with the
representing measure $\nu$ on $\partial H$ having the density
$1_{B(0,1)}(z)(1+|z|^2)^{-d/2} |z|^{-d} f(z)$. Notice that
$\nu(\{0\})=0$. Further, $f(0)=0$ (otherwise the integral in
(\ref{c:criterion2})
would diverge). Let $x\in A\cap B(0,1/2)$. Then
$x_d\le f(\wt{x},0)\le a|\wt{x}|$, and so $x_d/(2a) \le |\wt{x}|/2
\le 1/4$. Further, for $z\in B((\wt{x},0), \frac{x_d}{2a})\cap
\partial H$ we have
\begin{align*}
|z| \, &\le \,  |z-\wt{x}|+|\wt{x}| \, \le \, \frac{x_d}{2a}+|\wt{x}| \, \le \,
\frac{|\wt{x}|}{2}+|\wt{x}| \, = \, \frac32 |\wt{x}|<1\, ,\\
|x-z| \, &= \, (|\wt{x}-z|^2+x_d^2)^{1/2} \, \le \, \left(
\left(\frac{x_d}{2a}\right)^2+x_d^2\right)^{1/2} \, = \,
x_d(1+(2a)^{-2})^{1/2}\, ,\\
f(z) \, &> \, f(\wt{x},0)-\frac{x_d}{2} \, \ge \,  x_d-\frac{x_d}{2}
\, = \, \frac{x_d}{2}\, .
\end{align*}
Therefore, with $
\omega_{d-1}
$ denoting the volume of the
$(d-1)$-dimensional unit ball and using Theorem
\ref{t:martin-estimate}, we get
\begin{eqnarray*}
s(x)&\ge &C_5(4)^{-1}\, V(\delta_H(x))\int_{B((\wt{x},0),
\frac{x_d}{2a})\cap \partial H}|x-z|^{-d} |z|^{-d} f(z)\,
d\lambda'(z)\\
&\ge & C_5(4)^{-1}\, V(\delta_H(x)) \,
\frac{1}{x_d^d(1+(2a)^{-2})^{d/2}}\, \frac{1}{\left(\frac32
|\wt{x}|\right)^d}\,
\frac{x_d}{2}
\omega_{d-1} \left(\frac{x_d}{2a}\right)^{d-1}\\
&\ge & c_1 C_5(4)^{-1}
\omega_{d-1}\, V(\delta_H(x)) x_d^{-d}
|\wt{x}|^{-d} x_d x_d^{d-1}\\
&\ge & c_1 C_5(4)^{-1}
\omega_{d-1}\, \frac{V(\delta_H(x))}{|x|^{d}}\,\ge \, c_1 C_5(4)^{-2}
\omega_{d-1} M^H(x,0)\, ,
\end{eqnarray*}
where $c_1$ depends on $d$ and the Lipschitz constant $a$. This
proves that for every $x\in A\cap B(0,1/2)$ it holds that
$$
\frac{s(x)}{M^H (x,0)} \, \ge \,  c_1 C_5(4)^{-2}
\omega_{d-1} \,  > \,  \nu(\{0\})\, .
$$
Hence, by Proposition \ref{p:new}, $A$ is
minimally thin in $H$ with respect to
$X$ at 0. \qed

\bigskip

{\bf Acknowledgements.} We thank the referee for very helpful comments
on the first version of this paper.
\vspace{.1in}
\begin{singlespace}

\end{singlespace}

\end{doublespace}

\bigskip

{\bf Panki Kim}

Department of Mathematical Sciences and Research Institute of Mathematics,

Seoul National University,
San56-1 Shinrim-dong Kwanak-gu,
Seoul 151-747, Republic of Korea

E-mail: \texttt{pkim@snu.ac.kr}

\bigskip

{\bf Renming Song}

Department of Mathematics, University of Illinois, Urbana, IL 61801, USA

E-mail: \texttt{rsong@math.uiuc.edu}

\bigskip

{\bf Zoran Vondra{\v{c}}ek}

Department of Mathematics, University of Zagreb, Zagreb, Croatia

Email: \texttt{vondra@math.hr}

\end{document}